\documentclass{amsproc}
\usepackage{aspmproc}

\renewcommand\qedsymbol{$\blacksquare$}

\usepackage{amssymb}
\usepackage{amsmath}
\usepackage{cite}
\usepackage{mathrsfs}

\makeatletter

\address{\bigskip\hfil\begin{tabular}{l@{}}
            School of Mathematics and Statistics F07\\  
            University of Sydney,
            Sydney N.S.W. 2006.\hfill\qquad
                        {\tt mathas@maths.usyd.edu.au}\\  
            Australia.\hfill\qquad
                    {\tt www.maths.usyd.edu.au/u/mathas/}
          \end{tabular}}

\renewcommand{\labelenumi}{(\roman{enumi})}

\def\And{\text{\ and\ }}

\def\For{\text{\ for\ }}
\def\ForAll{\text{\ for all\ }}

\def\ForSome{\text{\ for some\ }}
\def\If{\text{\ if\ }}

\def\Whenever{\text{\ whenever\ }}

\def\Otherwise{\text{\ otherwise}}
\def\th{{\text{th}}}

\def\Set[#1]#2|#3|{\Big\{\ #2\ \Big| \
            \vcenter{\hsize #1mm\centering#3}\Big\}}

{\catcode`\|=\active
  \gdef\set#1{\mathinner{\lbrace\,{\mathcode`\|"8000%
                                   \let|\midvert #1}\,\rbrace}}
}
\def\midvert{\egroup\mid\bgroup}

\def\Number#1{\refstepcounter{equation}
              \leqno(\theequation)\if*#1%
              \else\def\@currentlabel{{\rm\theequation}}\label{#1}%
              \fi}
\def\Dag{\ifmmode\leqno(\dag)\else$(\dag\)$\fi}
\def\DDag{\ifmmode\leqno(\ddag)\else$(\ddag\)$\fi}

\numberwithin{equation}{section}

\newtheorem{Definition}[equation]{Definition}
\newtheorem{Theorem}[equation]{Theorem}
\newtheorem{Proposition}[equation]{Proposition}

\newtheorem{Corollary}[equation]{Corollary}
\newtheorem{Conjecture}[equation]{Conjecture}

\newtheorem{Problem}[equation]{Problem}

\let\SR\stackrel
\def\Prod{\displaystyle\prod}

\def\){\big)}
\def\({\big(}
\let\>\rangle
\let\<\langle

\let\iso\cong

\let\ss\subseteq
\let\0\varnothing              

\let\bar\overline

\let\gedom\trianglerighteq
\let\gdom\vartriangleright

\def\End{\mathop{\rm End}\nolimits}
\def\Hom{\mathop{\rm Hom}\nolimits}
\def\Ind{\mathop{\rm Ind}\nolimits}
\def\Res{\mathop{\rm Res}\nolimits}

\def\N{{\mathbb N}}
\def\Z{{\mathbb Z}}
\def\Q{{\mathbb Q}}

\def\C{{\mathbb C}}

\let\To\longrightarrow
\def\map#1#2{\,{:}\,#1\!\longrightarrow\!#2}
\def\mapsto{\!\longmapsto\!}
\def\bij#1#2{\,{:}\,#1\!\xrightarrow\sim\!#2}

\title[Representations of cyclotomic algebras]
{The representation theory of the Ariki-Koike and cyclotomic
$q$-Schur algebras}
\author{Andrew Mathas}


\usepackage{pstricks}
\psset{linewidth=0.8pt}  

\def\diag{\operatorname{diag}}

\def\rtuple#1{({#1}^{(1)},\dots,{#1}^{(r)})}

\def\1{\mathbf{i}}
\def\A{\mathcal A}
\def\B{\mathfrak B}
\def\D{\mathcal D}
\def\F{\mathcal F}
\def\H{\mathscr H}
\def\Haff{\hat H}
\def\L{\mathscr L}
\def\q{\mathbf Q}
\def\Sym{\mathfrak S}
\def\Schur{\mathscr S}

\def\G{\mathbf G}
\def\Bo{\mathbf B}

\def\Ind{\operatorname{Ind}}
\def\iInd{i\text{-}\Ind}
\def\Res{\operatorname{Res}}
\def\iRes{i\text{-}\Res}
\def\Soc{\operatorname{Soc}}

\let\phi\varphi

\def\[{[\![}
\def\]{]\!]}

\def\nr{\underline{\mathbf{nr}}}

\def\GL{\operatorname{GL}}

\def\mod{\text{-{\bf mod}}}
\def\proj{\text{-{\bf proj}}}

\let\len\ell
\def\comp{\operatorname{comp}}
\def\res{\operatorname{res}}
\def\Shape(#1){\operatorname{Shape}(#1)}
\def\Type(#1){\operatorname{Type}(#1)}

\def\Sfun{\Phi_\omega}

\def\p{\mathfrak p}
\def\O{\mathcal O}

\def\llen(#1){\ell\ell(#1)}

\def\Rad{\operatorname{Rad}}
\def\rest{{\downarrow}}
\def\s{\mathfrak s}
\def\t{\mathfrak t}

\let\uacc=\u
\def\u{\mathfrak u}
\def\v{\mathfrak v}

\def\a{\mathbf a}
\def\g{\mathfrak g}
\def\gl{\mathfrak{gl}}
\def\usl{U(\widehat{\mathfrak{sl}}_e)}
\def\uvsl{U_v(\widehat{\mathfrak{sl}}_e)}

\let\Sect\S
\def\S{\mathsf S}
\def\T{\mathsf T}
\def\U{\mathsf U}
\def\V{\mathsf V}
\def\sA{\mathsf A}
\def\sB{\mathsf B}

\def\phiST{\phi_{\S\T}}
\def\phio{\phi_\omega}

\def\Std(#1){{\mathcal T}^{\text{s}}\!(#1)}
\def\SStd(#1,#2){{\mathcal T}^{\text{ss}}_{#2}\!(#1)}

\def\myarrow#1{\,\ooalign{\raise.8ex\hbox{\ \small#1}\crcr%
           $\xrightarrow{\phantom{\hbox{\ \small#1}}}$}\,}
\def\iarrow{\myarrow{\ $i$}}
\let\arrow\longrightarrow

\def\goodarrow{\myarrow{\rm\ good}}

\usepackage[vcentermath]{youngtab}
\def\bitab(#1|#2){\bigg(\ \young(#1\ ,\quad\young(#2)\ \bigg)}
\def\tritab(#1|#2|#3){{\bigg(\begin{array}[t]{*3c}
                \begin{array}[t]{@{}l@{}}\young(#1)\end{array}\ ,&
                \begin{array}[t]{@{}l@{}}\young(#2)\end{array}\ ,&
                \begin{array}[t]{@{}l@{}}\young(#3)\end{array}
              \end{array}\bigg)}}

\newdimen\hoogte    \hoogte=10pt    
\newdimen\breedte   \breedte=14pt   
\newdimen\dikte     \dikte=0.7pt    

\newenvironment{Young}{\begingroup
       \def\vr{\vrule height0.8\hoogte width\dikte depth 0.2\hoogte}
       \def\fbox##1{\vbox{\offinterlineskip
                    \hrule height\dikte
                    \hbox to \breedte{\vr\hfill$##1$\hfill\vr}
                    \hrule height\dikte}}
       \vtop\bgroup \offinterlineskip \tabskip=-\dikte \lineskip=-\dikte
            \halign\bgroup &\fbox{##\unskip}\unskip  \crcr 
     }
     {\egroup\egroup\endgroup}

\def\Tritab(#1|#2|#3){\bigg(\ \begin{Young}#1\cr\end{Young}\,,\ 
                              \begin{Young}#2\cr\end{Young}\,,\ 
                              \begin{Young}#3\cr\end{Young}\ \bigg)}
\makeatother

\begin{document}
\bibliographystyle{andrew}
\maketitle

\section{Introduction}

The Ariki-Koike algebras were introduced by Ariki and
Koike~\cite{AK1} who were interested in them because they are a
natural generalization of the Iwahori-Hecke algebras of types $A$ and
$B$. At almost the same time, Brou\'e and Malle~\cite{BM:cyc} attached
to each complex reflection group a {\it cyclotomic Hecke algebra}
which, they conjectured, should play a role in the decomposition of
the induced cuspidal representations of the finite groups of Lie type.
The Ariki-Koike algebras are a special case of Brou\'e and Malle's
construction.

The deepest conjectures of Brou\'e and Malle (and others) concerning
the Ariki-Koike algebras have not yet been proved; however, many of
the consequences of these conjectures have been established.
Further, the representation theory of these algebras is beginning
to be well understood. For example, the simple modules of the
Ariki-Koike algebras have been classified; the blocks are known;
there are analogues of Kleshchev's modular branching rules; and, in
principle, the decomposition matrices of the Ariki-Koike algebras are
known in characteristic zero. In many respects this theory looks much
like that of the symmetric groups; in particular, there is a rich
combinatorial mosaic underpinning these results which involves
familiar objects like standard tableaux (indexed by multipartitions),
Specht modules and so on. 

The cyclotomic Schur algebras were introduced by Dipper, James and the
author~\cite{DJM:cyc}; by definition these algebras are endomorphism
algebras of a direct sum of ``permutation modules'' for an
Ariki-Koike algebra. This is generalizes the Dipper-James
definition~\cite{DJ:Schur} of the $q$-Schur algebras as endomorphism
algebras of tensor space. We were interested in these algebras both as
another tool for studying the Ariki-Koike algebras and because we
hoped that there might be a cyclotomic analogue of the famous
Dipper-James theory~\cite{DJ:Schur,BDK}. 

As with the Ariki-Koike algebras, the representation theory of the
cyclotomic Schur algebras is now well developed. They are always
cellular algebras; indeed, they are quasi-hereditary. The cellular
basis of these algebras is indexed by a generalization of semistandard
tableaux and their representation theory looks very much like the
representation theory of the $q$-Schur algebras. In particular, they
have a highest weight theory; there is a cyclotomic Schur functor and
a double centralizer theorem; the Jantzen filtrations of the
cyclotomic Weyl modules satisfy a generalization of the Jantzen sum
formula; and the cyclotomic Schur algebras have Borel subalgebras and
admit a triangular decomposition.

In the short time since its inception this theory has blossomed
producing many interesting results; largely this is because it
generalizes the representation theories of the symmetric groups, the
Schur algebras and the $q$-analogues of these. Many of the results in
this article have the flavour of results from Lie theory; however, as
yet, there are no known connections between the representation
theories of the cyclotomic Schur algebras and the finite groups of Lie
type except in the case where the underlying complex reflection group
is actually a Weyl group.

The aim of this article is to describe the representation theory of
these algebras in detail. Throughout we have tried to give an
indication of how the results are proved; unfortunately, in distilling
one or more papers in to one or more paragraphs some of the finer
details have been lost. 

\section{The Ariki-Koike algebras}

In this chapter we introduce the Ariki-Koike algebras by giving three
different constructions of them. From the point of view of
presentations it is clear that all three definitions agree; however,
for motivation, and also for proving certain results, it is important
to know the different contexts in which the Ariki-Koike algebras
arise. We remark that the Ariki-Koike algebras also appeared in a
different guise in the work of Cherednik~\cite{Cherednik:GT}.

We begin with a brief discussion of the complex reflection groups
which underpin the Ariki-Koike algebras. In the final section we give
a brief account of the conjectures of Brou\'e and Malle~\cite{BM:cyc}
which describe the role that the Ariki-Koike algebras should play in
the representation theory of the finite groups of Lie type.

\subsection{The complex reflection group of type $G(r,1,n)$.} Fix
integers $r\ge1$ and $n\ge0$ and let $W_{r,n}=\Z/r\Z\wr\Sym_n$ be the
wreath product of a cyclic group of order $r$ and a symmetric group of
degree $n$. Then $W_{r,n}$ is the complex reflection group of type
$G(r,1,n)$ in the Shephard-Todd classification~\cite{ST}; in
particular, $W_{r,n}$ has a faithful representation on a complex
vector space on which it acts as a group generated by reflections (see
section~(\ref{braid})).

If $r=1$ then $W_{1,n}\cong\Sym_n$ is just the symmetric group $\Sym_n$.
If $r=2$ then $W_{2,n}=\Z/2\Z\rtimes\Sym_n$ is the hyperoctohedral group,
or the group of signed permutations. In these two cases $W_{r,n}$ is a
{\it Coxeter group} or {\it real reflection group}; in fact, they are the
Weyl groups of type $A_{n-1}$ and $B_n$ respectively.

The group $W_{r,n}$ has the Coxeter like presentation given by
the following diagram.
$$\begin{array}{c}
\psline[doubleline=true,doublesep=.6mm](-1.75,0)(-1,0)
\psline(-1,0)(.7,0)
\psline(1.3,0)(2,0)
\cput(-2,0){r}
\qdisk(-1,0){4pt}
\qdisk(0,0){4pt}
\rput(1,0){\cdots}
\qdisk(2,0){4pt}
\rput(-2,-0.4){t_0}
\rput(-1,-0.4){t_1}
\rput(0,-0.4){t_2}
\rput(2,-0.4){t_{n-1}}\\[10pt]
\end{array}$$
The circle around the $r$ indicates that the corresponding generator
$t_0$ has order $r$; otherwise, this should be read as a standard Dynkin
diagram. Thus, as an abstract group, $W_{r,n}$ is generated by
elements $t_0,t_1,\dots,t_{n-1}$ which are subject to the relations
$$\begin{array}{r@{\,=\,}l@{\qquad}l}
t_0^r&1,                         \\
t_i^2&1,                         &\For 1\le i<n,\\
t_0t_1t_0t_1&t_1t_0t_1t_0,\\
t_it_j&t_jt_i,                   &\For 0\le j<i-1<n-1,\\
t_it_{i+1}t_i&t_{i+1}t_it_{i+1}, &\For 1\le i<n-1.
\end{array}$$
In particular, the subgroup $\<t_1,\dots,t_{n-1}\>$ of $W_{r,n}$
is isomorphic to the symmetric group~$\Sym_n$; hereafter, we identify
$\Sym_n$ and $\<t_1,\dots,t_{n-1}\>$ via the map $(i,i+1)\mapsto t_i$,
for $1\le i<n$.

Let $l_1=t_0$, $l_2=t_1t_0t_1,\dots, l_n=t_{n-1}\dots t_1t_0t_1\dots
t_{n-1}$. Then $l_1,\dots,l_n$ generate a subgroup of $W_{r,n}$
isomorphic to $\Z/r\Z\times\dots\times\Z/r\Z$  ($n$ copies), which is
just the base group when we consider $W_{r,n}$ as the semidirect
product $(\Z/r\Z\times\dots\times\Z/r\Z)\rtimes\Sym_n$. Thus, as a set,
$W_{r,n}=\set{l_1^{a_1}\dots l_n^{a_n} w| 0\le a_i<r\And w\in\Sym_n}$
and these elements are all distinct. In particular, $|W_{r,n}|=r^nn!$.

In general, $W_{r,n}$ is not a Coxeter group so the familiar
combinatorics of root systems and length functions cannot be used in
understanding $W_{r,n}$ and its representations. (Bremke and
Malle~\cite{BM:redwds} have defined a root system for $W_{r,n}$.) The
theory of complex reflection groups is still very much in its infancy;
the major tool being used to understand these groups is the geometry
of their reflection representation.

\subsection{The Ariki-Koike algebras} The Iwahori-Hecke algebras of
Weyl groups play an important role in the representation theory of the
groups of Lie type. Two important special cases of these algebras are
the Iwahori-Hecke algebras of the Weyl groups of types $A_{n-1}$
and~$B_n$ which are the groups $G(1,1,n)$ and $G(2,1,n)$,
respectively. Ariki and Koike~\cite{AK1} observed that the
definition of these algebras could be generalized to give a Hecke
algebra, or deformation algebra, for each complex reflection group of
type $G(r,1,n)$. 

Let $R$ be an integral domain with $1$ and let $q,Q_1,\dots,Q_r$ be
elements of $R$ with $q$ invertible. Let $\q=\{Q_1,\dots,Q_r\}$.

Deforming the relations of $W_{r,n}$ we obtain the Ariki-Koike
algebra.

\begin{Definition}
[Ariki-Koike~\cite{AK1}] 
The {\sf Ariki-Koike algebra}is the unital
associative $R$-algebra  $\H_{q,\q}(W_{r,n})$ with generators 
$T_0,T_1,\dots,T_{n-1}$ and relations
$$\begin{array}{r@{\,=\,}l@{\qquad}l}
(T_0-Q_1)\dots(T_0-Q_r)&1,       \\ 
(T_i-q)(T_i+1)&1,                &\For 1\le i<n,\\
T_0T_1T_0T_1&T_1T_0T_1T_0,\\
T_iT_j&T_jT_i,                   &\For 0\le i<j-1<n-1,\\
T_iT_{i+1}T_i&T_{i+1}T_iT_{i+1}, &\For 1\le i<n-1.
\end{array}$$
\end{Definition}
The three homogeneous relations are known as {\it braid relations}.

Typically, we write $\H=\H_{q,\q}(W_{r,n})$; when we wish to emphasize
the ring of definition we will write $\H=\H_{R,q,\q}(W_{r,n})$.

Notice that if $R$ contains a primitive $r^\th$ root of unity $\zeta$
and we set $q=1$ and $Q_s=\zeta^s$, for $1\le s\le r$, then 
$\H\cong R W_{r,n}$, the group algebra of $W_{r,n}$ (for this choice
of parameters the relations collapse to give those of $W_{r,n}$).

Let $w\in\Sym_n$. Then $w=t_{i_1}\dots t_{i_k}$ for some $i_j$ with
$1\le i_j<n$. If $k$ is minimal we say that $t_{i_1}\dots t_{i_k}$ is
a {\sf reduced expression} for $w$ and define $T_w=T_{i_1}\dots T_{i_k}$.
Since the braid relations hold in $\H$ it follows from Matsumoto's
monoid lemma (see, for example, \cite[Theorem~1.8]{M:ULect}), that
$T_w$ is independent of the choice of reduced expression for $w$.

Mimicking the definition of the elements $l_k$ in $W_{r,n}$, for
$k=1,\dots,n$ set $L_k=q^{1-k}T_{k-1}\dots T_1T_0T_1\dots T_{k-1}$.
(The renormalization by the unit $q^{1-k}$ is there to make the
combinatorics more natural later on.) Using the relations it is
straightforward to see that $L_1,\dots,L_n$ generate an abelian
subalgebra of $\H$ and that the symmetric polynomials in
$L_1,\dots,L_n$ belong to the centre of~$\H$.

{\it A priori} there is no reason to expect that the presentation
above will yield an interesting algebra. The first indication that
$\H$ is worth studying is the following theorem.

\begin{Theorem}
[Ariki-Koike~\cite{AK1}]
The Ariki-Koike algebra $\H$ is free as an $R$-module with
basis
$\set{L_1^{a_1}\dots L_n^{a_n} T_w| 0\le a_i<r\And w\in\Sym_n}$.
\label{AK basis}\end{Theorem}

In particular, notice that $\H$ is $R$-free of rank $r^nn!=|W_{r,n}|$
for any choice of $R$, $q$ and~$\q$. Furthermore, the subalgebra of
$\H$ generated by $T_1,\dots,T_{n-1}$ is isomorphic to the
Iwahori-Hecke algebra $\H_q(\Sym_n)$ of the symmetric
group $\Sym_n$. Hereafter, we identify the two algebras $\H_q(\Sym_n)$
and $\<T_1,\dots,T_{n-1}\>$.

Using the relations it is not hard to show that  $\H$ is spanned by
the elements $L_1^{a_1}\dots L_n^{a_n} T_w$; there are $|W_{r,n}|$
such elements. To prove linear independence Ariki and Koike explicitly
constructed the simple $\H$-modules using a generalization of Young's
seminormal form for the Ariki-Koike algebras when
$R=\C(q,u_1,\ldots,u_r)$; see Theorem~\ref{AK} below.  This shows that
$\H/\Rad\H$ has dimension at most $|W_{r,n}|$. Hence, $\H$ is
semisimple and Theorem~\ref{AK basis} is proved when
$R=\C(q,Q_1,\dots,Q_r)$. The general case now follows by a
specialization argument.

There are now other proofs of Theorem~\ref{AK basis} available.  Brou\'e,
Malle and Rouquier~\cite[Theorem~4.24]{BMR} have given a geometrical
argument which results from thinking of $\H$ as a quotient of the
group algebra of the braid group of $W_{r,n}$ and studying its
monodromy representation; this is the topic of the next section.
Sakamoto and Shoji~\cite{SakSho} also proved Theorem~\ref{AK basis} as a
consequence of an analogue of Schur-Weyl reciprocity for $\H$ and a
particular quantum group; we will return to this in section~\ref{SS}
below. 

Finally, we remark that Shoji~\cite{Sho:Frob} has given a different
presentation of $\H$ when $Q_1,\dots,Q_r$ are all
distinct (for $R$ an integral domain).

\subsection{The braid group of $W_{r,n}$ and the Hecke algebra}
\label{braid}

At almost the same time that Ariki and Koike introduced their algebra,
Brou\'e and Malle~\cite{BM:cyc} associated to each complex reflection
group~$W$ a {\it cyclotomic Hecke algebra}; for the group $W_{r,n}$
Brou\'e and Malle's cyclotomic Hecke algebra is precisely the
Ariki-Koike algebra. Brou\'e and Malle's motivation was that they
expected that the cyclotomic Hecke algebras should play a role in the
representation theory of the finite groups of Lie type similar to, but
more complicated than, that played by the Iwahori-Hecke algebras
(see section~2.5).

In this section we briefly describe Brou\'e and Malle's definition in
the case of $W_{r,n}$ and some of its consequences.

Let $V$ be the complex vector space with basis
$\{e_1,\dots,e_n\}$ and let $\zeta\in\C$ be a primitive $r^\th$ root
of unity. The symmetric group $\Sym_n=\<t_1,\dots,t_{n-1}\>$ acts on
$V$ in the natural way; extend this to an action of $W_{r,n}$ by
letting $t_0$ act via the $n\times n$ matrix $\diag(\zeta,1,\dots,1)$.
This defines a faithful representation of $W_{r,n}$. Observe that each
of the generators of $W_{r,n}$ acts as a {\it reflection} (that is,
fixes a space of codimension~$1$), so this shows that $W_{r,n}$ is a
complex reflection group. 

Let $\Omega=\set{e_i-\zeta^ke_j|1\le j\le i\le n\And\max(j-i,-1)<k<r}$.
Then $\Omega$ is in one-to-one correspondence with the set of
reflections in $W_{r,n}$, the correspondence attaching to each reflection
an eigenvector with non-trivial eigenvalue;
see~\cite[\Sect3]{BM:redwds}. For each~$\omega\in\Omega$ let $H_\omega$
be the hyperplane orthogonal to $\omega$ and let $\mathscr M=V\setminus
\bigcup_{\omega\in\Omega} H_\omega$ be the associated hyperplane
complement and $\mathscr M/W_{r,n}$ its quotient by $W_{r,n}$. 

\begin{Definition}
The {\sf braid group} of $W_{r,n}$ is the group
$$\B_{r,n}=\pi_1(\mathscr M/W_{r,n},x_0),$$
where $x_0\in\mathscr M/W_{r,n}$.
\end{Definition}

Here, $\pi_1(\mathscr M/W_{r,n},x_0)$ is the fundamental group of the
quotient space $\mathscr M/W_{r,n}$ with base point~$x_0$. Because
$\mathscr M$ is connected $\B_{r,n}$ is independent of the choice of
$x_0$.

If $r>1$ then $\B_{r,n}$ is a braid group of type $B_n$ and as an
abstract group it is generated by elements $s_0,\dots,s_{n-1}$ subject
to the relations
$$s_0s_1s_0s_1=s_1s_0s_1s_0,\quad
s_is_j=s_js_i,\quad\And\quad
s_is_{i+1}s_i=s_{i+1}s_is_{i+1},
$$
where $1\le i<n-1$, $0\le j<n$ and $|i-j|>1$. In particular, observe
that~$W_{r,n}$ is a quotient of $\B_{r,n}$ (via the map which sends $s_i$
to $t_i$ for $0\le i<n$).

The generators of $\B_{n,r}$ can be chosen as generators of the
monodromy around the hyperplanes. Bessis~\cite{Bessis} has now given a
general argument for the existence of such presentations for the braid
groups of complex reflection groups. 

Brou\'e and Malle considered the algebra $R\B_{r,n}/I_{q,\q}$,
where~$I_{q,\q}$ is the ideal of $R\B_{r,n}$ generated by
$(s_0-Q_1)\dots(s_0-Q_r)$ and $(s_i-q)(s_i+1)$, for $1\le i<n$;
evidently, $\H\cong R\B_{r,n}/I_{q,\q}$.  One consequence of this
definition is that we can use the monodromy representation
of the braid group $\B_{r,n}$ to analyze $\H$. This leads to a more
conceptual proof of the fact that $\H$ is always free as an
$R$-module of rank $|W_{r,n}|$ (a corollary of Theorem~\ref{AK basis}).
Moreover, it yields the following important result.

\begin{Theorem}
[\protect{Brou\'e-Malle-Rouquier~\cite[Theorem 4.24]{BMR}}] 
Let $\mathbb K=\C(q,Q_1,\dots,Q_r\)$. Then the monodromy
representation of $\B_{r,n}$ induces an isomorphism of 
$\mathbb K$-algebras $\H_{\mathbb K,q,\q}\cong\mathbb KW_{r,n}$. 
\label{iso}\end{Theorem}

Here, $\mathbb KW_{r,n}$ is the group algebra of $W_{r,n}$
over $\mathbb K$. That $\H_{\mathbb K,q,\q}$ and $\mathbb KW_{r,n}$
are isomorphic algebras can be established by a general Tits
deformation theory argument (see, for example, \cite[\Sect66]{C&R:2});
an explicit isomorphism is also given in
\cite[Corollary~3.13]{M:gendeg}. The main point of this result is that
the isomorphism is canonically determined.

Lusztig~\cite{L:iso} has proved a similar result for the
Iwahori-Hecke algebras of Weyl groups; however, his argument is less
elementary relying on a deep property of the cells of Weyl groups. For
Weyl groups, Lusztig's isomorphism and that of Theorem~\ref{iso} are
different.

\subsection{The affine Hecke algebra of type $A$}\label{affine H} The
Ariki-Koike algebras should really be considered as affine objects
because they are quotients of the (extended) affine Hecke algebra of
type $A$ (i.e., the affine Hecke algebra of $\GL_n(\C)$). The affine
Hecke algebra $\Haff_n$ is the $R$-algebra with generators
$T_1,\dots,T_{n-1}$ and
$X_1^{\pm1},\dots,X_n^{\pm1}$ and relations
$$\begin{array}{ccc}
(T_i-q)(T_i+1)=1,& T_iT_{i+1}T_i=T_{i+1}T_iT_{i+1}, 
                 & T_iX_iT_i=qX_{i+1}\\
T_iT_k=T_kT_i,   & X_iX_k=X_kX_i,  & T_iX_k=X_kT_i
\end{array}$$ 
and $X_iX_i^{-1}=1=X_i^{-1}X_i$ for all sensible values of $i,j,k$
with $|i-k|>1$.  In particular, abusing notation slightly, notice that
there is surjective algebra homomorphism $\Haff_n\twoheadrightarrow\H$
given by sending $T_i\mapsto T_i$ and $X_j\mapsto L_j$, for $1\le i<n$
and $1\le j\le n$ respectively. It is easy to see that
$\H_{q,\q}(W_{r,n})\cong\Haff_n/\<(X_1-Q_1)\dots(X_1-Q_r)\>$.

It follows from the relations that $T_1,\dots,T_{n-1}$
generate a subalgebra of $\Haff_n$ isomorphic to $\H_q(\Sym_n)$ and that
$X_1^{\pm1},\dots,X_n^{\pm1}$ generate a Laurent polynomial ring. 
Therefore, as an $R$-module,  
$\Haff_n\cong\H_q(\Sym_n)\otimes R[X_1^{\pm1},\dots,X_n^{\pm1}]$;
so, $\Haff_n$ is a twisted tensor product.

Let $\mathcal P=\bigoplus_{i=1}^n\Z\epsilon_i$ be the free
$\Z$-module with basis $\epsilon_1,\dots,\epsilon_n$; so, 
$\mathcal P$ is the weight lattice of $\GL_n(\C)$. The symmetric group
$\Sym_n$ acts on $\mathcal P$ by permuting the $\epsilon_i$. 

If $\lambda\in\mathcal P$ set $X^\lambda=X_1^{\lambda_1}\dots
X_n^{\lambda_n}$. Then the two commutation relations for the $T_i$
and the~$X_j$ can be replaced by the relation
$$ T_i X^\lambda=X^{t_i\lambda}T_i
              +(q-1)\frac{X^\lambda-X^{t_i\lambda}}{1-X^{\alpha_i}},$$
where $\lambda\in\mathcal P$, $\alpha_i=\epsilon_i-\epsilon_{i+1}$ and
$1\le i<n$. A quick calculation shows that $X^\lambda-X^{t_i\lambda}$ is
divisible by $1-X^{\alpha_i}$ so the right hand side does make sense.
Notice that when $q=1$ this relation becomes 
$t_i X^\lambda=X^{t_i\lambda}t_i$; this is what we expect because
the extended affine Weyl group is the semidirect product $\mathcal
P\rtimes\Sym_n$.

Now suppose that $R$ is an algebraically closed field. Bernstein
showed that the centre of~$\Haff_n$ is the set of symmetric
polynomials in $X_1,\dots,X_n$ (Theorem~\ref{Bernstein}).
Consequently, $\Haff_n$ is finite dimensional over its centre;
therefore, by Schur's Lemma, every irreducible $\Haff_n$-module is
finite dimensional (with dimension at most $n!$ since
$\dim_R\Haff_n/Z(\Haff_n)=(n!)^2$ by Theorem~\ref{Bernstein} below). 

As remarked above, each Ariki-Koike algebra $\H_{q,\q}(W_{r,n})$ is a
quotient of $\Haff_n$, so every irreducible $\H$-module is also an
irreducible $\Haff_n$-module.  Conversely, if $M$ is an irreducible
$\Haff_n$-module let $c_M(X_1)$ be the characteristic polynomial for the
action of~$X_1$ on~$M$; then $\H_M:=\Haff_n/\<c_M(X_1)\>$ is an
Ariki-Koike algebra (with parameters the eigenvalues for the action of
$X_1$ on~$M$) and $M$ is an irreducible $\H_M$-module. (More
generally,~$M$ is an irreducible module for any Ariki-Koike algebra
obtained by quotienting out by the ideal generated by any polynomial in
$X_1$ which is divisible by $c_M(X_1)$.) Thus the irreducible
$\Haff_n$-modules are precisely the irreducible
$\H_{q,\q}(W_{r,n})$-modules as $\q$ ranges over the elements of
$(R^\times)^r$ for~$r\ge1$.

\subsection{The conjectures of Brou\'e, Malle and Michel} The
conjectures which we now discuss grew out of the attempts of Brou\'e
and others to understand Brou\'e's~\cite{Broue:perfect} conjectures
for blocks with abelian defect groups in the case of the finite
reductive groups. We consider only a very special case of these
conjectures; for references and further details see the original
papers \cite{BM:cyc,BM:cohom,BMM} and
Brou\'e's~\cite{Broue:conjectures} comprehensive survey article.

Let $\G$ be an algebraic group defined over $\bar{\mathbb F}_q$, where
$q$ is a prime power, and let $W$ be the Weyl group of $\G$. Let
$F\map\G\G$ be a Frobenius map and let $G=\G^F$ be the $F$-fixed
points of $\G$. We assume that $W$ is $F$-split. The simplest example
is to take $\G=\GL_n(\bar{\mathbb F}_q)$ and $F(a_{ij})=(a_{ij}^q)$;
then $G=\GL_n(q)$ and $W=\Sym_n$.

Let $\Bo$ be an $F$-stable Borel subgroup of $\G$ and set $B=\Bo^F$.
It is well-known that the irreducible constituents of $\Ind_B^G(1)$
are in one-to-one correspondence with the irreducible
representations of $W$; see, for example, \cite{Carter:2}. The
Iwahori-Hecke algebras of Weyl groups play an important role in this
theory; indeed, $\H_q(W)\cong\End_G\(\Ind_B^G(1)\)$ and this explains
why the dimensions of the irreducible representations in the unipotent
principal series, the constituents of $\Ind_B^G(1)$, are given by
evaluating certain polynomials $D_\chi(x)$ at $x=q$. The conjectures
which follow attempt to explain other ``generic'' features of the
representation theory of finite reductive groups.

Let $B_W$ be the Braid group of $W$ and for $w\in W$ let $\underline
w\in B_W$ be the lift of $w$ (under the canonical embedding of $W$
into the positive braid monoid $B_W^+$). Brieskorn and
Saito~\cite{BrieSaito} showed that the centre of $B_W$ is generated by
$\pi=\underline{w_0^2}$ (or $\underline{w_0}$ if $w_0$ is central in
$W$), where $w_0$ is the unique element of maximal length in $W$. 

Call an element $w\in W$ {\sf good} if $\pi=\underline{w}^d$ for some
$d$. Note that $w$ has order $d$ since $w_0^2=1$ in $W$. Every
conjugacy class of regular elements in~$W$ contains a good element.
Assume that $w$ is good. Then Brou\'e and Michel~\cite{BM:cohom}
have shown that every good element is regular; so $C_W(w)$ is a complex
reflection group by Springer~\cite{Springer:reg}. Let $B_w=B(C_W(w))$
be the braid group of $C_W(w)$. It is conjectured that
$B_w=C_{B_W}(\underline w)$; this has now been proved in almost all
cases \cite{BDM,Bessis:dual}.

Let $X_w$ be the Deligne-Lusztig variety associated to $w$; so $X_w$
is the variety of Borel subgroups $\Bo$ of $\G$ such that $\Bo$ and
$F(\Bo)$ are in relative position $w$. Fix a prime $\ell$ not dividing
$q$ and consider the \'etale cohomology groups
$H^i_c(X_w,\bar\Q_\ell)$ of $X_w$. The finite group $G=\G^F$ acts on
$X_w$ and hence also on $H^i_c(X_w,\bar\Q_\ell)$. By
\cite{BM:cohom,DMR} there is also an action of $C_{B_W}(\underline w)$
on $H^i_c(X_w,\bar\Q_\ell)$ (this comes from an action of the positive
braids in $C_{B_W}(\underline w)$ on $X_w$).  In many cases the action
of $\bar\Q_\ell C_{B_W}(\underline w)$ is known to factor through a
cyclotomic Hecke algebra.  Conjecturally, the action of
$C_{B_W}(\underline w)$ should generate 
$\End_{\bar\Q_\ell G}\(H^i_c(X_w,\bar\Q_\ell)\)$; this seems to be
very hard.

Let $\mathcal H(\G,F,W,w)$ be the image of $\bar\Q_\ell B_w$ in the
(graded) endomorphism algebra of 
$\bigoplus_{i\ge0} H^i_c(X_w,\bar\Q_\ell)$. Then  
$\mathcal H(\G,F,W,w)$ is a finite dimensional algebra and the
following should be true.

\begin{Conjecture}
[Brou\'e,Malle,Michel~\cite{Broue:conjectures,BM:cyc,BM:cohom,BMM}]
\quad\newline\label{Broue's conjectures}
Suppose that $w$ is a good element of order $d$.
\begin{enumerate}
\item If $i\ne j$ then the $\bar\Q_\ell\G^F$-modules
$H^i_c(X_w,\bar\Q_\ell)$ and $H^j_c(X_w,\bar\Q_\ell)$ have no
irreducible constituents in common.
\item There is a $d$-cyclotomic Hecke algebra
$\H_x\(C_W(w)\)$ of the complex reflection group $C_W(w)$ such that
$$\mathcal H(\G,F,W,w)\cong\H_{\bar\Q_\ell,q}\(C_W(w)\)\cong
 \End_{\bar\Q_\ell\G^F}\Big(\bigoplus_{i\ge0}H^i_c(X_w,\bar\Q_\ell)\Big).
$$
\item There is a one-to-one correspondence
$\chi\longleftrightarrow\chi_q$ between the irreducible
representations of $C_W(w)$ and the irreducible constituents
of the $\bar\Q_\ell\G^F$-module $\bigoplus_{i\ge0}H^i_c(X_w,\bar\Q_\ell)$. Moreover, for each
irreducible character $\chi$ of $C_W(w)$ there is a polynomial
$D_\chi(x)$, depending only on $\chi$,  such that the degree of
$\chi_q$ is equal to $D_\chi(q)$.
\end{enumerate}\end{Conjecture}

We now explain the term {\it $d$-cyclotomic Hecke algebra} when
$C_W(w)=W_{r,n}$. Let $\xi\in\C$ be a root of unity and $x$ an
indeterminate over $\Z[\xi]$ and let $\zeta_d$ be a primitive $d^\th$
root of unity. An Ariki-Koike algebra
$\H_{R,v,\q}(W_{r,n})=\H_x(W_{r,n})$ is {\it $d$-cyclotomic} if
$R=\Z[\xi][x,x^{-1}]$ and it has parameters of the form
$v=\zeta^{a_v}x^{b_v}$ and $Q_s=\zeta^{a_s}x^{b_s}$, for some rational
numbers $a_v,b_v$ and $a_s,b_s$, such that:
{\renewcommand{\labelenumi}{(\alph{enumi})}
\begin{enumerate}
\item $\H_{\zeta_d}=\H_x(W_{r,n})\otimes_R R/(x-\zeta_d) 
                      \cong\Z[\xi]W_{r,n}$; and, 
\item $\H_q=\H_x(W_{r,n})\otimes_RR/(x-q)$ is semisimple
over its field of fractions.  
\end{enumerate}}\noindent%
(For example, take parameters $v=x^d$ and $Q_s=x^{s-1}$ (with $\xi=1$);
then $\H_{\zeta_d}\cong\Z[\zeta_d]W_{r,n}$.)

Thus, part (ii) of the conjecture together with (b) implies that the
irreducible representations occurring in $\bigoplus_i
H^i_c(X_w,\bar\Q_\ell)$ are in one-to-one correspondence with the
irreducible representations of $\H_q\(C_W(w)\)$; in turn, by (a) these
representations are in one-to-one correspondence with the irreducible
representations of $C_W(w)$. Importantly, nothing here depends upon
the choice of~$q$ or~$\ell$.  Conjecturally, these correspondences
come from a derived equivalence, so they are really perfect isometries
(``bijection with signs'').  The polynomials $D_\chi(x)$ in part
(iii) are the generic degrees of $\H_x\(C_W(w)\)$; see the remarks
after Theorem~\ref{Schur elements}.

In fact, part (iii) of the conjecture is known to be true. The main
reason why this is known is that the virtual module
$\bigoplus_{i\ge0}(-1)^iH^i_c(X_w,\bar\Q_\ell)$ is a Deligne-Lusztig
representation (specifically, it is $R_{T_w}^G(1)$, where $T_w$ is the
maximal torus associated to the conjugacy class of $w$ in $W$), so its
irreducible constituents are known. Parts (i) and (ii) of the
conjecture are known in only a small number of cases.

We also mention that everything above is compatible with the
decomposition of the unipotent characters of $\G^F$ into
$d$-Harish-Chandra series~\cite{BMM}. For these details, and
stronger forms of the conjecture, we refer the reader to Brou\'e's
article~\cite{Broue:conjectures}.

To conclude this section we remark that if $w=1$ then $X_1=G/B$ is
the flag variety; so, $H_c^0(X_1,\bar\Q_\ell)\cong\Ind_B^G(1)$ and
all higher cohomology groups are zero. Thus, in this case the
conjectures recover the well-known results for the principal
unipotent series of~$\G^F$. According to our definitions, $w=1$ is
not a good element of $W$; however, we have discussed only a special
case of the general conjectures. 

\section{The representation theory of the Ariki-Koike algebras}

\subsection{The semisimple representation theory of $\H$}

Because $W_{r,n}$ is the wreath product $\Z/r\Z\wr\Sym_n$, its
ordinary irreducible representations are indexed by $r$-tuples of
partitions of~$n$. In this section we see that the same is true of the
irreducible representations of $\H$ when $\H$ is semisimple.

A {\sf partition} of $n$ is a sequence
$\sigma=(\sigma_1\ge\sigma_2\ge\cdots)$ of non-negative integers
$\sigma_i$ such that $|\sigma|=\sum_{i\ge1}\sigma_i=n$; we write
$\sigma=(\sigma_1,\dots,\sigma_k)$ if $\sigma_i=0$ for $i>k$. A
{\sf multipartition} of~$n$ is an ordered $r$-tuple
$\lambda=\rtuple\lambda$ of partitions with
$|\lambda^{(1)}|+\dots+|\lambda^{(r)}|=n$. We write $\lambda\vdash n$
if $\lambda$ is a multipartition of $n$.

The set of multipartitions is a poset under {\sf dominance} $\gedom$,
where $\lambda\gedom\mu$ if
$$\sum_{t=1}^{s-1}|\lambda^{(t)}|+\sum_{j=1}^i\lambda^{(s)}_j
   \ge\sum_{t=1}^{s-1}|\mu^{(t)}|+\sum_{j=1}^i\mu^{(s)}_j$$
for $s=1,2,\dots,r$ and all $i\ge1$. If $\lambda\gedom\mu$ and
$\lambda\ne\mu$ we write $\lambda\gdom\mu$.

The {\sf diagram} of $\lambda$ is 
$[\lambda]=\set{(i,j,s)|1\le j\le\lambda^{(s)}_i\And 1\le s\le r}.$
The elements of $[\lambda]$ are called {\sf nodes}; more generally, a
node is any triple $(i,j,s)$ where $1\le s\le r$ and~$i,j\ge1$.

A $\lambda$-tableau is a bijection $\t\map{[\lambda]}\{1,2\dots,n\}$,
which we consider as an $r$-tuple $\t=\rtuple\t$ of labeled tableaux
where $\t^{(s)}$ is a $\lambda^{(s)}$-tableau for each $s$; the
tableaux $\t^{(s)}$ are the {\sf components} of $\t$. If $\t$ is a
$\lambda$-tableau we write $\Shape(\t)=\lambda$. 

A tableau $\t$ is {\sf standard} if, in each component, its entries
increase from left to right along each row and from top to bottom in
down each column. For example,
$$\t^\lambda=\tritab(123,4|5,6|78,9) \And
  \t=\tritab(479,6|1,8|25,3)\Number{tab example}$$
are two standard $\((3,1),(1^2),(2,1)\)$-tableaux. Let
$\Std(\lambda)$ be the set of standard $\lambda$-tableaux.

If $\t$ is a $\lambda$-tableau and $w\in\Sym_n$ let $\t w=\t\circ w$
be the tableau obtained from $\t$ by replacing each entry in $\t$ by
its image under $w$. This defines a right action of $\Sym_n$ on the
set of all $\lambda$-tableaux. For example,
$\t=\t^\lambda(1,4,6,8,5)(2,7)(3,9)$ in (\ref{tab example}).

If $\t$ is a tableau and $k$ an integer, with $1\le k\le n$, then the
{\sf residue} of $k$ in $\t$ is defined to be $\res_\t(k)=q^{j-i}Q_s$,
where $k$ appears in row $i$ and column $j$ of $\t^{(s)}$; that is,
$\t(i,j,s)=k$.

The last ingredient that we need is something like the Poincar\'e
polynomial of a Coxeter group; however, be warned that it is not true
that $|W_{r,n}|=P_\H(q,\q)$ when $R=\C$, $q=1$ and~$Q_s=\zeta^s$
(that is, when $\H=RW_{r,n}$). Let
$$P_\H(q,\q)=\prod_{i=1}^n (1+q+\dots+q^{i-1})\cdot \prod_{1\le
i<j\le r}\prod_{-n<d<n}(q^dQ_i-Q_j).$$ 

We can now describe the irreducible representations of $\H$ when
$P_\H(q,\q)$ is invertible. (Note that if $R$ is a field then
$P_\H(q,\q)$ is invertible if and only if $P_\H(q,\q)\ne0$.)

\begin{Theorem}
[Ariki-Koike~\cite{AK1}]\label{H ss reps}\label{AK}
Suppose that $P_\H(q,\q)$ is invertible in $R$.
\begin{enumerate}
\item For each multipartition $\lambda$ let $V^\lambda$ be the
$R$-module with basis 
$$\set{v_\t|\t\text{\ a standard $\lambda$-tableau}}.$$ 
Then $V^\lambda$ becomes a right $\H$-module via 
$v_\t T_0=\res_\t(1)v_\t$ and, for $1\le i<n$, if $\s=\t t_i$ is not 
standard then
\begin{align*}
v_\t T_i&=\begin{cases} 
  \ qv_\t,&\text{if $i$ and $i+1$ are in the same row of $\t$,}\\
  -v_\t,&\text{if $i$ and $i+1$ are in the same column of $\t$,}\\
\end{cases}
\intertext{and if $\s$ is standard then}
v_\t T_i&=
  \frac{(q-1)\res_\t(i)}{\res_\t(i)-\res_\s(i)}v_\t
  +\frac{q\res_\t(i)-\res_\s(i)}{\res_\t(i)-\res_\s(i)}v_\s.
\end{align*}
\item If $R$ is a field then $V^\lambda$ is an irreducible 
$\H$-module for each multipartition $\lambda$.
\item If $R$ is a field then $\set{V^\lambda|\lambda\vdash n}$ is a 
complete set of pairwise non-isomorphic irreducible $\H$-modules.
\end{enumerate}\end{Theorem}

Part (i) is proved by a brute force calculation to show that the
action of the generators on~$V^\lambda$ respects the relations in
$\H$. The remaining parts can be proved by looking at how the
commutative subalgebra $\L=\<L_1,\dots,L_n\>$ of $\H$ acts on
$V^\lambda$. From Theorem~\ref{AK}(i) it follows that $v_\t L_k=\res_\t(k)v_\t$
for all standard tableaux $\t$, for $1\le k\le n$. Hence $R v_\t$ is
an irreducible $\L$-module; in fact, Ariki and
Koike~\cite{AK1} show that every irreducible $\L$-module is
of this form. Moreover, because $P_\H(q,\q)\ne0$ if~$\s$ and $\t$ are
standard tableaux then $\s=\t$ if and only if $\res_\t(k)=\res_\s(k)$,
for $1\le k\le n$; this implies that $Rv_\s\cong Rv_\t$ as
$\L$-modules if and only if~$\s=\t$. As a consequence,
$S^\lambda\cong S^\mu$ if and only if $\lambda=\mu$. Part~(iii) now
follows by counting dimensions because
$$\dim\H\ge\dim(\H/\Rad\H)\ge\sum_\lambda(\dim V^\lambda)^2
   =r^nn!=|W_{r,n}|\ge\dim\H.$$
(The third equality follows from the Robinson-Schensted
correspondence which implies that the sum of the squares of the number
of standard $\lambda$-tableaux, as $\lambda$ runs over all
multipartitions of $n$,  is equal to~$|W_{r,n}|$.) As we have equality
throughout, this also proves Theorem~\ref{AK basis} (indeed, this is
how Ariki and Koike first proved it).

\begin{Corollary}
[Ariki~\cite{Ariki:ss}] Suppose that $R$ is a field. Then 
$\H$ is semisimple if and only if $P_\H(q,\q)\ne0$.
\label{H ss}
\end{Corollary}

\begin{proof}
[Sketch of proof]
By Theorem~\ref{AK} if $P_\H(q,\q)\ne0$ then $\H$ is semisimple. For
the converse, when $P_\H(q,\q)=0$ the ideal of~$\H$ generated by
$$\Big(\prod_{k=1}^n\prod_{s=1}^{r-1}(L_k-Q_s)\Big)
      \Big(\sum_{w\in\Sym_n}T_w\Big),$$
is nilpotent. (This ideal affords the ``trivial'' representation 
of~$\H$.)
\end{proof}

Halverson and Ram~\cite{HalRam} have generalized the
Murnaghan-Nakayama rule of the symmetric groups to give a method for
computing the characters of the irreducible
representations~$V^\lambda$.  (In fact, they also compute the
characters of the irreducible representations of the cyclotomic Hecke
algebras of type $G(r,p,n)$; the irreducible representations of these
algebras were constructed by Ariki~\cite{Ariki:hecke}.) See also
Shoji~\cite{Sho:Frob}.

As remarked earlier the symmetric polynomials in $L_1,\dots,L_n$
belong to the centre of~$\H$. In the semisimple case this is a
complete description of the centre.

\begin{Theorem}
[Ariki-Koike~\cite{AK1}] Suppose that $R$ is a field and that
$P_\H(q,\q)\ne0$. Then the centre of $\H$ is equal to the set of
symmetric polynomials in $L_1,\dots,L_n$.  
\label{H symmetric}\end{Theorem}

Graham~\cite{Graham:centre} has recently shown that the centre of
$\H_{R,q}(\Sym_n)$ is always equal to the set of symmetric polynomials
in $L_1,\dots,L_n$ when $R$ is an integral domain (this is the case
$r=1$). Ariki~\cite{Ariki:can} has given an example which shows that
the centre of $\H$ can be larger than the set of symmetric polynomials
when $r>1$.

When $P_\H(q,\q)\ne0$ the author~\cite{M:gendeg} has explicitly
described the primitive central idempotents as symmetric polynomials
in $L_1,\dots,L_n$ (see also Shoji~\cite{Sho:Frob}); this gives a
second proof of Theorem~\ref{H symmetric}. In addition,
\cite{M:gendeg} constructs the primitive idempotents and a Wedderburn
basis of $\H$ in the semisimple case.

Define $\tau\map\H R$ to be the $R$-linear map determined by
$$\tau(L_1^{a_1}\dots L_n^{a_n} T_w)
        =\begin{cases} 1,&\If a_1=\dots=a_n=0\And w=1,\\
                0,&\Otherwise,
\end{cases}$$
for $0\le a_i<r$ and $w\in\Sym_n$.  Notice that if $q=1$ and
$Q_s=\zeta^s$, where $\zeta=\exp(2\pi i/r)\in\C$, then $\tau$ is the
natural trace function on the group algebra $\C W_{r,n}$. The
definition of $\tau$ looks quite arbitrary; however, as we explain
below, $\tau$ is canonical.

\begin{Proposition}
\label{symmetric}
Assume that $R$ is an integral domain. Then the following hold.
\begin{enumerate}
\item$(${\bf Bremke-Malle}~\cite{BM:redwds}$)$ $\tau$ is a trace 
form on $\H$.
\item$(${\bf Malle-Mathas}~\cite{MM:trace}$)$ Suppose that
$q,Q_1,\dots,Q_r$ are all invertible in $R$. Then $\tau$ is 
non-degenerate. Consequently, $\H$ is a symmetric algebra.
\end{enumerate}\end{Proposition}

Part (i) is straightforward; although we should mention that Bremke
and Malle use a different (but, by \cite{MM:trace}, equivalent),
definition of the trace form $\tau$. For the Iwahori-Hecke algebras
($r\le 2$), part (ii) is also routine (see, for example,
\cite[Prop.~1.16]{M:ULect}); in contrast, whilst not difficult, the
proof of~(ii) is a laborious calculation when $r>2$. As an indication
of the difficulties here, no pair of dual bases for $\H$ is known when
$r>2$ (except in the semisimple case;
see~\cite[Theorem~3.9]{M:gendeg}).

As we will describe, Proposition~\ref{symmetric} provides the
strongest known link between the representation theory of $\H$ and
that of the finite groups of Lie type (when $r>2$).

If $R$ is a field and $P_\H(q,\q)\ne0$ then $\H$ is semisimple. Let
$\chi^\lambda$ be the character of $V^\lambda$.  Since $\tau$ is a
trace function we can write 
$$\tau=\sum_{\lambda\vdash n}\frac1{\s_\lambda(q,\q)}\chi^\lambda$$
for some $s_\lambda(q,\q)\in R$. The rational functions
$s_\lambda(q,\q)$ are the {\sf Schur elements} of $\H$; to describe
them we need some more notation.

Define the length of a partition $\sigma$ to be the smallest integer
$\len(\sigma)$  such that $\sigma_i=0$ for all~$i>\len(\sigma)$;
the {\sf length} of a multipartition $\lambda$ is 
$\len(\lambda)=\max\set{\len(\lambda^{(s)})|1\le s\le r}$.
Suppose that $L\ge\len(\lambda)$ and set 
$\beta^{(s)}_i=\lambda^{(s)}_i+L-i$ for $i=1,\dots,L$ and 
$1\le s\le r$; also set $B_s=\{\beta^{(s)}_1,\dots,\beta^{(s)}_L\}$ 
for $s=1,\dots,r$. The matrix $B=\(\beta^{(s)}_i)_{s,i}$ is the
(ordinary) $L$-{\sf symbol} of~$\lambda$~\cite{BK, Malle:gendeg}. 

\begin{Theorem}
[Geck-Iancu-Malle~\cite{GIM}] 
Suppose that $\lambda$ is a multipartition of $n$ with
$L$-symbol $B=(\beta^{(s)}_i)_{s,i}$ such that $L\ge\len(\lambda)$. 
Then $s_\lambda(q,\q)$ is equal to 
$$(-1)^{a_{rL}}q^{b_{rL}}
\frac{\Prod_{1\le s<t\le r}(Q_s-Q_t)^L\ \cdot
      \prod_{1\le s,t\le r}\prod_{\alpha_s\in B_s}
      \prod_{1\le k\le\alpha_s} (q^kQ_s-Q_t)}
     {(q-1)^n(Q_1\dots Q_r)^n
      \Prod_{1\le s\le t\le r}
      \prod_{\substack{(\alpha_s,\alpha_t)\in B_s\times B_t\\
                      \alpha_s>\alpha_t\If s=t}}
             (q^{\alpha_s}Q_s-q^{\alpha_t}Q_t)},$$
where $a_{rL}=n(r-1)+\binom r2\binom L2$ 
and $b_{rL}=\frac{rL(L-1)(2rL-r-3)}{12}$.
\label{Schur elements}\end{Theorem}

For $r=1,2$ the Schur elements were first computed by
Hoefsmit~\cite{Hoefsmit}. Murphy~\cite{murphy:hecke} gave a different
argument for type $A$ (that is, $r=1$). For $r>2$ this result was
conjectured by Malle~\cite{Malle:gendeg}. Geck, Iancu and Malle use a
clever specialization argument due to Orellana~\cite{Orel:B} to
compute the Schur elements using the Markov trace of the Hecke
algebras $\H_q(\Sym_m)$; in turn, this builds on work of
Wenzl~\cite{Wenzl:subfactors}. It is not hard to see that if
$f_\lambda$ is a primitive idempotent in $\H$ which generates the
Specht module $S^\lambda$ then $s_\lambda(q,\q)=\tau(f_\lambda)$; this
observation is used in~\cite{M:gendeg} to give a direct proof of 
Theorem~\ref{Schur elements}.

Theorem~\ref{Schur elements} is important because when combined with
\cite[3.16 and~6.11]{Malle:gendeg} it implies that
$\Phi_d$-blocks~\cite{BMM} of the finite reductive groups satisfy a
generalized Howlett--Lehrer theory~\cite{HL:hecke}; more precisely,
Conjecture~\ref{Broue's conjectures}(iii) is true and the
dimensions of the irreducible representations in a unipotent
$\Phi_d$-block are given by specializations of the generic
degrees of $\H$; these are the polynomials
$D_\lambda(q)=s_\eta(q,\q)/s_\lambda(q,\q)$, where
$\eta=((n),(0),\dots,(0))$.

As a second application of Theorem~\ref{Schur elements}, it follows
from \cite[Theorem~5.2]{GIM} and \cite[Lemma~2.7]{BMM:spetses} that
the trace form $\tau$ is the unique trace form on $\H$ which, in a
precise sense~\cite[Theorem~2.1]{BMM:spetses}, is compatible with the
usual trace forms on both $W_{r,n}$ and on the braid group $\B_{r,n}$.
In addition, Malle~\cite{Malle:expgendeg} uses Theorem~\ref{Schur
elements} to define the notion of ``spetsiality'' for complex
reflection groups; for more details see \cite{BMM:spetses}. 

Finally, Brou\'e and Kim~\cite{BK} use Theorem~\ref{Schur elements},
together with the block structure of~$\H$, to show that the
irreducible representations of $\H$ can be grouped according to a
generalization of Lusztig's {\it famillies}; a key ingredient in their
paper is a block theoretical characterisation of Lusztig's famillies
due to Rouquier~\cite{Rouq:blocks}. Again, the combinatorial
description of the spetsial famillies of~$\H$ had previously been
conjectured by Malle~\cite{Malle:gendeg}.

\subsection{The modular representation theory of $\H$} We now turn to
the modular representation theory of $\H$; that is, the representation
theory when $\H$ is not semisimple. In types $A$ and $B$ the
irreducible modular representations were first constructed by Dipper
and James~\cite{DJ:reps} and Dipper, James and Murphy~\cite{DJM},
respectively; Graham and Lehrer~\cite{GL} considered the general case
using cellular algebra techniques. Even though the papers
\cite{murphy:basis,DJM} predated Graham and Lehrer, the cellular
approach is already implicit in them.

Graham and Lehrer constructed a cellular basis for $\H$ by building
upon the Kazhdan-Lusztig basis of $\H_q(\Sym_n)$ (which is itself
cellular).  We will describe a different cellular basis of $\H$ which
comes from the work of Dipper, James and the author~\cite{DJM:cyc}.
We prefer this basis because we know how to lift this basis to give a
basis for the cyclotomic $q$-Schur algebras and because this basis
has many nice combinatorial and representation theoretic properties.

Fix a multipartition $\lambda$ and let 
$\Sym_\lambda
    =\Sym_{\lambda^{(1)}}\times\dots\times\Sym_{\lambda^{(r)}}$ 
be the associated Young subgroup of $\Sym_n$. Equivalently,
$\Sym_\lambda$ is the row stabilizer of the $\lambda$-tableau
$\t^\lambda$ which has the numbers $1,\dots,n$ entered in order from
left to right, top to bottom first along the rows
of~$\t^{\lambda^{(1)}}$ and then $\t^{\lambda^{(2)}}$ and so on (for
example, see the first tableau in (\ref{tab example})).

Define elements $x_\lambda$ and $u_\lambda^+$ in $\H$ by
$$x_\lambda=\sum_{w\in \Sym_\lambda} T_w\quad\And\quad
  u_\lambda^+=\prod_{s=2}^r
      \prod_{k=1}^{|\lambda^{(1)}|+\dots+|\lambda^{(s-1)}|}(L_k-Q_s).$$
It follows easily from the relations in $\H$ that 
$x_\lambda u_\lambda^+=u_\lambda^+ x_\lambda$. Set 
$m_\lambda=x_\lambda u_\lambda^+$. Although somewhat ungainly, the
function of $u_\lambda^+$ is used to control the eigenvalues of the $L_k$
on the modules below.

If $\t$ is a standard $\lambda$-tableau let $d(\t)\in\Sym_n$ be the
unique permutation in $\Sym_n$ such that $\t=\t^\lambda d(\t)$. 

Let $*$ be the anti-isomorphism of $\H$ determined by
$T_i^*=T_i$ for $0\le i<n$; then $*$ is an involution and
$T_w^*=T_{w^{-1}}$, $L_k^*=L_k$ and $(h_1h_2)^*=h_2^*h_1^*$.

\begin{Definition}
Suppose that $\s$ and $\t$ are standard $\lambda$-tableaux. Let
$m_{\s\t}=T_{d(\s)}^*m_\lambda T_{d(\t)}$.
\end{Definition}

\begin{Theorem}[The standard basis theorem~\cite{DJM:cyc}]
The Ariki-Koike algebra $\H$ is free as an $R$-module with cellular
basis
$$\set{m_{\s\t}|\s\And\t\text{\ standard $\lambda$-tableaux,
                      $\lambda\vdash n$}}.$$
\label{Murphy basis}\end{Theorem}

When $r=1$ this result is due to Murphy~\cite{murphy:basis} and when
$r=2$ it was proved by Dipper, James and Murphy~\cite{DJM}. The basis
$\{m_{\s\t}\}$ is called both the {\sf Murphy basis} and the {\sf standard
basis} of $\H$. As mentioned above, Graham and Lehrer~\cite{GL} were
the first to produce a cellular basis of $\H$.

The proof of this theorem starts by observing that $\H$ is spanned by
a set of more general elements $m_{\s\t}$ where $\s$ and $\t$ are {\it
row standard} tableaux of the same shape. (The entries in row standard
tableaux increase along rows, but not necessarily down columns.) Next,
one shows that if $\s$ and $\t$ are not standard tableaux then
$m_{\s\t}$ can be written as a linear combination of ``higher terms''
$m_{\u\v}$; so, by induction, $\H$ is spanned by standard basis
elements (here, ``higher'' is essentially the Bruhat order on
$\Sym_n$). The rewriting rules involved in this step are essentially
Garnir relations; in fact, they are a little bit easier than the
classical Garnir relations because we work modulo a filtration. A
counting argument now shows that we have a basis. In order to show
that the basis is cellular some accounting details need to be carried
through the argument; this adds only minor complications to the proof.

We will not describe the theory of cellular algebras in detail;
instead the reader is referred to the beautiful paper of Graham and
Lehrer~\cite{GL}. Two other approaches to cellular algebras can be
found in \cite{KX:CellAlg,M:ULect}.

The required indexing of a cellular basis is already implicit in our
notation. The two properties that the basis $\{m_{\s\t}\}$
must satisfy for it to be cellular are: (i) the $R$-linear map
determined by $m_{\s\t}\mapsto m_{\t\s}$ must be an algebra
anti-isomorphism --- this is obvious for us because
$m_{\s\t}^*=m_{\t\s}$; and, (ii) for all  $\lambda$-tableaux~$\t$ and
all $h\in\H$ there exist scalars $r_\v\in R$ such that for any standard
$\lambda$-tableau $\s$
$$m_{\s\t}h\equiv
  \sum_{\v\in\Std(\lambda)} 
   r_\v m_{\s\v}\pmod{\mathcal H^\lambda},\Number{cellular}$$ 
where $\mathcal H^\lambda$ is the $R$-module spanned by the elements
$m_{\u\v}$ for $\Shape(\u)=\Shape(\v)\gdom\lambda$. The point of this
equation is that the scalars $r_\v$ depend only on $\t,\v$ and $h$;
importantly, $r_\v$ does not depend on $\s$.

Applying the anti-isomorphism $*$ to the last equation gives a
left hand analogue of (\ref{cellular}) for $hm_{\s\t}$.
It follows that $\mathcal H^\lambda$ is a two-sided ideal of $\H$.

\begin{Definition}
Suppose that $\lambda$ is a multipartition of $n$. The {\sf Specht
module} $S^\lambda$ is the right $\H$-module generated by 
$m_\lambda+\mathcal H^\lambda$.
\label{Specht modules}\end{Definition}

Thus, $S^\lambda$ is a submodule of the quotient module 
$\H/\mathcal H^\lambda$. Du and Rui~\cite{DuRui:branching} have shown
how to construct the Specht modules as submodules of $\H$ (as distinct
from subquotients as we have defined them here).

For each standard $\lambda$-tableau $\t$ let
$m_\t=m_{\t^\lambda\t}+\mathcal H^\lambda
     =m_\lambda T_{d(\t)}+\mathcal H^\lambda$. 
It follows from Theorem~\ref{Murphy basis} that $S^\lambda$ 
is free as an $R$-module with basis 
$\set{m_\t|\t\text{\ a standard $\lambda$-tableau}}$; moreover,
by (\ref{cellular}) the action of $\H$ on this basis is given by
$$m_\t h=\sum_{\substack{\v\text{\ standard}\\\text{$\lambda$-tableau}}}
      r_\v m_\v,$$
where the scalars $r_\v\in R$ are the same as those in
(\ref{cellular}).  It follows from the left and right handed versions
of (\ref{cellular}) that there is a bilinear form on~$S^\lambda$ which is
determined by 
$$\<m_\s,m_\t\>m_{\u\v}\equiv m_{\u\s}m_{\t\v}\pmod{\mathcal H^\lambda}$$
for all standard $\lambda$-tableaux~$\s$ and $\t$. This form is
associative in the sense that $\<xh,y\>=\<x,yh^*\>$ for all~$x,y\in
S^\lambda$ and $h\in\H$. Hence, $\Rad S^\lambda=\set{x\in
S^\lambda|\<x,y\>=0\ForAll y\in S^\lambda}$ is a submodule
of~$S^\lambda$ and we may make the following definition.

\begin{Definition}
Suppose that $\lambda$ is a multipartition of $n$. Then
$D^\lambda$ is the right $\H$-module 
$D^\lambda=S^\lambda/\Rad S^\lambda$.
\end{Definition}

Everything that we have said since Theorem~\ref{Murphy basis} is part of the
general machinery of cellular algebras. Without too much work, the
cellular theory now produces the following result.

\begin{Theorem}
[Graham-Lehrer~\cite{GL}, Dipper-James-Mathas~\cite{DJM:cyc}]
\label{main H theorem}
Suppose that $R$ is a field.
\begin{enumerate}
\item For each multipartition $\mu$,
$D^\mu$ is either zero or absolutely irreducible.
\item $\set{D^\mu|\mu\vdash n\And D^\mu\ne0}$ is a 
complete set of pairwise non-isomorphic irreducible $\H$-modules.
\item If $D^\mu\ne0$ then the decomposition multiplicity
$[S^\lambda{:}D^\mu]\ne0$ only if $\lambda\gedom\mu$; further,
$[S^\mu{:}D^\mu]=1$.
\end{enumerate}\end{Theorem}

Graham and Lehrer proved this result for a different collection of
modules; but this should really be considered their result. Again, for
the cases $r=1,2$ see \cite{DJ:reps,DJM}.

In particular, note that every field is a splitting field for $\H$.
The reader might be concerned with the claim that any field is a
splitting field for $\H$ because, for example,~$\Q$ is not a splitting
field for $W_{r,n}$ when $r>2$; however, this is OK because by
definition all of the eigenvalues of $T_0$ automatically belong
to~$R$.

The multiplicities $d_{\lambda\mu}=[S^\lambda{:}D^\mu]$ are the 
{\sf decomposition numbers} of~$\H$ and the matrix $(d_{\lambda\mu})$
is the decomposition matrix of $\H$.  Part~(iii) of Theorem~\ref{main
H theorem} says that the decomposition matrix of $\H$ is unitriangular
when its rows and columns are ordered in a way that is compatible with
the dominance order.

Corollary~\ref{H ss} and the theory of cellular algebras also gives us the
following result. 

\begin{Theorem}
Suppose that $R$ is a field. Then the following are equivalent.
\begin{enumerate}
\item $P_\H(q,\q)\ne0$;
\item $\H$ is semisimple;
\item $\H$ is split semisimple; and,
\item $S^\lambda=D^\lambda$ for all multipartitions $\lambda$ of $n$.
\end{enumerate}\end{Theorem}

If $1\le k\le n$ let $\t\rest k$ be the subtableau of~$\t$ containing
the integers $1,2,\dots,k$; so, if $\t$ is standard then
$\Shape(\t\rest k)$ is a multipartition of $k$.  We extend the
dominance ordering to the set of standard tableaux by defining
$\s\gedom\t$ if $\Shape(\s\rest k)\gedom\Shape(\t\rest k)$ for
$k=1,\dots,n$. Again we write $\s\gdom\t$ if~$\s\gedom\t$ and
$\s\ne\t$.

A useful fact about the standard basis of $\H$ is the following.

\begin{Proposition}
[\protect{\cite[Prop.~3.7]{JM:cyc-Schaper}}]\label{L action}
Suppose that $1\le k\le n$ and let $\s$ and $\t$ be standard tableaux
of the same shape. Then, there exist scalars $r_\v\in R$ such that
$$m_{\s\t}L_k=\res_\t(k)m_{\s\t}+\sum_{\v\gdom\t}r_\v m_{\s\v}
        \pmod{\mathcal H^\lambda}.$$
\end{Proposition}

As shown in \cite{JM:cyc-Schaper}, the general case can be reduced to
the case $r=1$ where it is a theorem of Dipper and
James~\cite{DJ:blocks}. When $r=1$ the result can be proved by
induction on~$n$ and $k$ using the fact that $L_1+\dots+L_n$ belongs
to the centre of $\H$; see~\cite{M:ULect}. 

As an application of Proposition~\ref{L action}, if $R$ is a field and
$P_\H(q,\q)\ne0$ then we can construct the irreducible $\H$-modules
either as the modules $V^\lambda$ of Theorem~\ref{H ss reps} or as the
Specht modules $S^\lambda$.  By Proposition~\ref{L action} the
modules~$V^\lambda$ and $S^\lambda$ have the same $\L$-module
composition factors; this implies that $V^\lambda\cong S^\lambda$ as
$\H$-modules.

We close this section with a reduction theorem which shows that, up to
Morita equivalence, the only important Ariki-Koike algebras are those
with parameters of the form (i) $Q_s=q^{a_s}$ for some integers $a_s$
with $|a_s|<n$, for $1\le s\le r$, or (ii) $Q_s=0$ for $1\le s\le r$.
The result actually says that we can reduce to the case where
there exists a constant $c\in R$ and integers $a_s$ such that
$Q_s=cq^{a_s}$, for all $s$; However, if $c\ne0$ then we can
renormalize the generator $T_0$ as $\tilde T_0=c^{-1}T_0$ and then the
order relation for $\tilde T_0$ becomes 
$(\tilde T_0-q^{a_1})\dots(\tilde T_0-q^{a_r})=0$, so we are back in 
case~(i).

Recall that $\q=\{Q_1,\dots,Q_r\}$ and fix a partition
$\q=\q_1\coprod\dots\coprod\q_\kappa$ (disjoint
union) of~$\q$ and let 
$$P_n(q,\q_1,\dots,\q_\kappa)=\prod_{1\le \alpha<\beta\le\kappa}
    \prod_{\substack{Q_i\in\q_\alpha\\ Q_j\in\q_\beta}}
    \prod_{-n<d<n}(q^dQ_i-Q_j).$$
Observe that $P_n(q,\q_1,\dots,\q_\kappa)$ is a factor of the
polynomial $P_\H(q,\q)$.

\begin{Theorem}
[Dipper-Mathas~\cite{DM:Morita}]\label{Morita}
Suppose that $R$ is an integral domain and that
$\q=\q_1\coprod\dots\coprod\q_\kappa$ is a partition of $\q$ such 
that the polynomial $P_n(q,\q_1,\dots,\q_\kappa)$ is invertible in $R$.
For $\alpha=1,\dots,\kappa$ let $r_\alpha=|\q_\alpha|$. Then
$\H_{q,\q}(W_{r,n})$ is Morita equivalent to the $R$-algebra 
$$\bigoplus_{\substack{n_1,\dots,n_\kappa\ge0\\n_1+\dots+n_\kappa=n}}
\H_{q,\q_1}(W_{r_1,n_1})\otimes\dots
          \otimes\H_{q,\q_\kappa}(W_{r_\kappa,n_\kappa}).$$
\end{Theorem}

If $r=2$ then $|\q_1|=|\q_2|=1$ and this is a result of Dipper and
James~\cite{DJ:B}. Du and Rui~\cite{DuRui:akmorita} extended the
argument of \cite{DJ:B} to prove the special case of Theorem~\ref{Morita} when
$|\q_\alpha|=1$ for $1\le\alpha\le\kappa$; notice that in this case
$\H$ is Morita equivalent to a direct sum of tensor products of
Iwahori-Hecke algebras of type $A$.

For the proof of Theorem~\ref{Morita} observe that by induction it is enough
to consider the special case $\kappa=2$. Without loss of generality we
may assume that $\q_1=\{Q_1,\dots,Q_s\}$ and
$\q_2=\{Q_{s+1},\dots,Q_r\}$ for some $s$. The trick is to consider
the right ideals~$V^b=v_b\H$, for $0\le b\le n$, where
$$v_b=\prod_{t=1}^s(L_1-Q_t)\dots(L_{n-b}-Q_t)\cdot T_{w_b}\cdot
      \prod_{t=s+1}^r(L_1-Q_t)\dots(L_b-Q_t)$$
and $w_b=(n,\dots,2,1)^b$. It turns out that the standard basis of
$\H$ can be adapted to give a `standard' basis of $V^b$. With this
basis in hand one sees that $V^b$ is a projective $\H$-module and
that $\End_\H(V^b)
\cong\H_{q,\q_1}(W_{s,b})\otimes\H_{q,\q_2}(W_{r-s,n-b})$ and
$\Hom_\H(V^b,V^c)=0$ for~$b\ne c$. These results imply that
$\bigoplus_{b=0}^n V^b$ is a projective generator for $\H$ which gives
the result. The Morita equivalence can be described very explicitly;
as one consequence, when $R$ is a field, it is easy to compare the
dimensions of the simple modules under the equivalence.

\subsection{Ariki's theorem}This section discusses a very deep result
of Ariki~\cite{Ariki:can} which gives a way to compute the
decomposition numbers of the Ariki-Koike algebras
$\H_{\C,q,\q}(W_{r,n})$ when $q\ne1$ and $Q_s\ne0$ for all $s$.
Throughout we assume that $R$ is a field.  It will be convenient to
write $\H_n=\H_{q,\q}(W_{r,n})$. Also let $\H_n$\mod\ be the category
of finite dimensional right $\H_n$-modules. We begin with some
motivation.

If $M$ is an $\H_n$-module let $\Res M$ be the
restriction of $M$ to $\H_{n-1}$. Then $\Res$ is an exact
functor from $\H_n$\mod\ to $\H_{n-1}$\mod. Since $\H_n$ is free
as an $\H_{n-1}$-module $\Res$ has a right adjoint; namely, the
induction functor which sends a right $\H_{n-1}$-module~$N$ to
$\Ind N=N\otimes_{\H_{n-1}}\H_n$.

If $\lambda$ is a multipartition of $n-1$ and $\mu$ is a
multipartition of $n$ write $\lambda\arrow\mu$ if the diagrams of
$\lambda$ and $\mu$ differ by only one node. From the definition of
the Specht modules it is clear that the action of $\H_{n-1}$ on $\Res
S^\mu$ is given by ignoring the node in the tableaux with label $n$.
With only a small amount of work this implies the following result.

\begin{Proposition}
[\protect{Ariki~\cite[Lemma~2.1]{Ariki:can}}]
Suppose that $\mu$ is a multipartition of $n$. Then $\Res S^\mu$ has a
filtration with composition factors isomorphic to the Specht
modules~$S^\lambda$, where $\lambda$ runs over the multipartitions
of $n-1$ such that $\lambda\arrow\mu$.
\label{branching}\end{Proposition}

Let $K_0(\H_n\mod)$ be the Grothendieck group of $\H_n\mod$. Thus,
$K_0(\H_n\mod)$ is the free abelian group generated by all finitely
generated right $\H_n$-modules where the relations are given by short
exact sequences. If $M$ is a right $\H_n$-module let $[M]$ be the
corresponding equivalence class in $K_0(\H_n\mod)$. By 
Theorem~\ref{main H theorem} $\set{[S^\mu]|D^\mu\ne0}$ and
$\set{[D^\mu]|D^\mu\ne0}$ are both bases of $K_0(\H_n\mod)$ and the
transition matrix between these bases is the decomposition matrix
of $\H$.

The functors $\Res$ and $\Ind$ induce homomorphisms of
abelian groups which, by abuse of notation, we also denote by $\Res$ 
and $\Ind$. Thus,
$\Res\map{K_0(\H_n\mod)}K_0(\H_{n-1}\mod)$ and
$\Ind\map{K_0(\H_n\mod)}K_0(\H_{n+1}\mod)$ are the maps
given by $\Res[M]=[\Res M]$ and $\Ind[M]=[\Ind M]$. These
homomorphisms are completely determined by their actions on the
Specht modules and this is given by Proposition~\ref{branching} 
and Frobenius reciprocity.

\begin{Corollary}
Suppose that $\lambda$ is a multipartition of $n$. Then 
$$\Res[S^\lambda]=\sum_{\nu\arrow\lambda}[S^\nu]
\quad\And\quad
\Ind[S^\lambda]=\sum_{\lambda\arrow\mu}[S^\mu].$$
\label{res}\end{Corollary}

Let $c_n=L_1+\dots+L_n$; then $c_n$ belongs to the centre of $\H_n$.
If $M$ is any $\H_n$-module let $M_\alpha=\set{m\in
M|(c_n-\alpha)^km=0\For k\gg0}$ be the corresponding generalized
eigenspace for $c_n$ acting on $M$, for $\alpha\in R$. Then $M_\alpha$
is an $\H_n$-module since $c_n\in Z(\H_n)$; so 
$M=\oplus_{\alpha\in R}M_\alpha$ as an $\H_n$-module.

Until further notice we assume that $q\ne1$ and that $Q_s=q^{a_s}$ for
some integers $a_s$, for $1\le s\le r$. In particular, this implies
that the eigenvalues of $c_n$ are always linear combinations of powers
of $q$. Let $e$ be the multiplicative order of $q$; then
$e\in\N\cup\{\infty\}$. 

Now the Specht module $S^\lambda$ is irreducible when $R=\C(q)$;
therefore, it follows from Proposition~\ref{L action}, and a
specialization argument, that $c_n$ acts on the Specht module
$S^\lambda$ as multiplication by the scalar
$c(\lambda)=\sum_{k=1}^n\res_{\t^\lambda}(k)$. Therefore,
$S^\lambda=(S^\lambda)_{c(\lambda)}$ is a single generalized
eigenspace and, by the Corollary,
$\Res S^\lambda=\oplus_{i\in\Z}(\Res S^\lambda)_{c(\lambda)-q^i}$
and
$\Ind S^\lambda=\oplus_{i\in\Z}(\Ind S^\lambda)_{c(\lambda)+q^i}$.
It follows that the eigenvalues of $c_n$ on an arbitrary $\H_n$-module
change by $\pm q^i$, for some $i\in\Z/e\Z$, under the functors $\Res$
and $\Ind$ respectively. Accordingly, we define new functors $\iRes$
and $\iInd$ on $\H_n\mod$ by
$$\iRes M=\bigoplus_\alpha (\Res M_\alpha)_{\alpha-q^i}\quad\And\quad
  \iInd M=\bigoplus_\alpha (\Ind M_\alpha)_{\alpha+q^i},
$$
for $i=0,1,\dots,e-1$. Then $\Res=\sum_{i=0}^{e-1}\iRes$ and
$\Ind=\sum_{i=0}^{e-1}\iInd$. These functors also induce group
homomorphisms $K_0(\H_n\mod)\To K_0(\H_{n\pm1}\mod)$ and these maps
are completely determined by their actions on the Specht modules.
Write $\lambda\iarrow\mu$ if $\lambda\arrow\mu$ and the node in
$[\mu]\setminus[\lambda]$ has residue $q^i$. Then we have the
following refinement of Corollary~\ref{res}.

\begin{Corollary}
Suppose that $0\le i<e$ and let 
$\lambda$ be a multipartition of $n$. Then 
$$\iRes[S^\lambda]=\sum_{\nu\iarrow\lambda}[S^\nu]
\quad\And\quad
\iInd[S^\lambda]=\sum_{\lambda\iarrow\mu}[S^\mu].$$
\end{Corollary}

Let $\H_n\proj$ be the category of finitely generated
{\it projective} $\H_n$-modules and let $K_0(\H_n\proj)$ be its
Grothendieck group. If $P$ is a projective
$\H_n$-module let $\[P\]$ denote its image in $K_0(\H_n\proj)$. 
Observe that there is a natural non-degenerate paring 
$$\<\ ,\ \>\map{K_0(\H_n\proj)\times K_0(\H_n\mod)}\Z$$
given by $\<\[P\],[M]\>=\dim_R\Hom_{\H_n}(P,M)$; hence,
$K_0(\H_n\proj)\cong K_0(\H_n\mod)^*$. Consequently, if $P^\mu$ is the
projective cover of~$D^\mu$ then $\set{\[P^\mu\]|\mu\vdash n\And
D^\mu\ne0}$ is a basis of~$K_0(\H_n\proj)$ and we have induced maps
$\iRes^*,\iInd^*\map{K_0(\H_n\proj)}K_0(\H_{n\pm1}\proj)$. 

We are almost ready to state Ariki's theorem. Let $\usl$ be the
Kac-Moody Lie algebra of type $A^{(1)}_{e-1}$. Thus, $\usl$ is the
$\C$-algebra generated by $d, e_i,f_i$ and $h_i$, for $0\le i<r$,
subject to a well-known set of relations; see \cite{Ariki:book,Kac}.
Let $\Lambda_0,\dots,\Lambda_{e-1}$ be the fundamental weights of
$\usl$ and recall that for each dominant weight
$\Lambda\in\sum_{i=0}^{e-1}\N\Lambda_i$ there is a unique integrable
highest weight $\usl$-module $L(\Lambda)$ with highest
weight~$\Lambda$.

\begin{Theorem}
[Ariki~\cite{Ariki:can,AM:simples}]\label{Ariki}
Suppose that $R$ is a field and fix $q,Q_1=q^{a_1},\dots,Q_r=q^{a_r}$
in~$R$ such that $q\ne1$ is a primitive $e^\th$ root of unity and
integers $a_1,\dots,a_r$ $($with $0\le a_i<e$ if $e<\infty)$. Finally,
let $\Lambda=\sum_{i=0}^{e-1}a_i\Lambda_i$ and set 
$V_{q,\q}(R)=\bigoplus_{n\ge0}K_0(\H_{R,n}\proj)\otimes_\Z\C$.
\begin{enumerate}
\item $V_{q,\q}(R)$ is an integrable $\usl$-module upon which the
Chevalley generators $e_i$ and $f_i$ act as follows:
$$\qquad\qquad
e_i\[M\]=\iRes^*\[M\]\quad\And\quad f_i\[M\]=\iInd^*\[M\],$$ for all 
$\[M\]\in V_{q,\q}(R)$. Moreover, $V_{q,\q}(R)\cong L(\Lambda)$ as a
$\usl$-module.
\item If $R$ is a field of characteristic zero then the canonical
basis of $V_{q,\q}(R)$ coincides with the basis
$$\set{\[P^\mu\]|D^\mu\ne0\ForSome \mu\vdash n\ge0}$$
given by the projective indecomposable 
$\H_n$-modules.
\end{enumerate}\end{Theorem}

Some remarks are in order. First, the hard part of this theorem is the
case where $R=\C$; this is proved in \cite{Ariki:can}. The result for
an arbitrary field follows from the complex case by a modular
reduction argument; see \cite{AM:simples}. Next, by the canonical
basis of $L(\Lambda)$ we mean the specialization at $v=1$ of the
Kashiwara-Lusztig canonical basis%
\footnote{Canonical bases of quantum groups were introduced independently
 by Lusztig~\cite{L:can} and Kashiwara~\cite{Kash:gcrys}.
 Jantzen~\cite{Jantzen:QG} has given an excellent treatment of this
 theory; unfortunately, he only considers quantum groups of finite type
 which is insufficient for our purposes.  Ariki~\cite{Ariki:book} gives a
 largely self-contained account of the canonical bases of $\uvsl$, which
 is exactly what we need. See also Lusztig's book~\cite{L:QGBK}.
}
of $L_v(\Lambda)$, the corresponding integrable highest weight
representation of the quantum group $\uvsl$. 

Theorem~\ref{Ariki} is a very deep result which relies upon the topological
$K$-theory of Kazhdan and Lusztig~\cite{KL:DelLang} and Ginzburg's
equivariant $K$-theory~\cite{ChrissGinz}; these theories give
different constructions of the standard modules of the affine Hecke
algebras in characteristic zero. For details of the proof see Ariki's
original paper~\cite{Ariki:can} and also his forthcoming
book~\cite{Ariki:book}. Geck~\cite{Geck:Ariki} has also written an
excellent survey article on the modular representation theory of Hecke
algebras; he includes a detailed account of Ariki's paper.

The special case of Theorem~\ref{Ariki} with $r=1$ proves the conjecture of
Lascoux, Leclerc and Thibon~\cite{LLT} for computing the decomposition
matrices of the Iwahori-Hecke algebras $\H_q(\Sym_n)$ of the
symmetric groups. The main point of \cite{LLT} is that they gave an
elementary combinatorial algorithm for computing the canonical basis
of the integrable highest weight module $L_v(\Lambda_0)$ for $\uvsl$
--- and hence the decomposition matrices of $\H_q(\Sym_n)$.  This
and similar algorithms are described in
\cite{Ariki:book,LLT,M:ULect,GW,LT:Schur2}. In contrast to the
difficulty of Theorem~\ref{Ariki}, these algorithms involve only
basic linear algebra; they amount to computing certain parabolic
affine Kazhdan-Lusztig polynomials of type $A$ and evaluating them at
$1$. This is described explicitly in \cite{LT:Schur2,GW}.

Uglov~\cite{Uglov},  extending the ideas of Leclerc and
Thibon~\cite{LT:Schur2}, has given an algorithm for computing the
canonical basis of any integrable highest weight module for $\uvsl$;
see also \cite{TakeUglov}. Hence, combining Theorem~\ref{Ariki}(ii) with
Uglov's work and Theorem~\ref{Morita} we have the following.

\begin{Corollary}
\label{decomp nos}
Suppose that $R$ is a field of characteristic zero and that
$q\ne1$ and $Q_s\ne0$ for $1\le s\le r$. Then the decomposition
matrix of $\H_{R,q,\q}(W_{r,n})$ is known.
\end{Corollary}

In practice there is a bit of work to be done to use this result to
compute the decomposition numbers of $\H$. First, Uglov's algorithm
computes a canonical basis for a larger space which contains
$L_v(\Lambda)$ as a submodule; this is less efficient than the LLT
algorithm and its variants. Next, Uglov's indexing of the canonical basis
of $L_v(\Lambda)$ is not compatible with Theorem~\ref{main H theorem}(ii) and
Theorem~\ref{H simples} below; a bijection between the different indexing sets
for the irreducibles is given by the paths in the associated crystal
graphs. Finally, the effect of the Morita equivalence of Theorem~\ref{Morita} on
the decomposition numbers must be taken into account; this last step is
straightforward and is described in \cite{DM:Morita}.

\subsection{The irreducible $\H$-modules} In principle, the simple
$\H_n$-modules are completely determined by Theorem~\ref{main H theorem};
that is, the simple $\H_n$-modules are precisely the non-zero
modules $D^\mu$ for $\mu$ a multipartition of~$n$.  Unfortunately, it
is non-trivial to determine when~$D^\mu$ is zero and when it is
non-zero.

We begin the classification of the simple modules of the Ariki-Koike
algebras with the case $r=1$.  Let $e$ be the smallest positive
integer such that $1+q+\dots+q^{e-1}=0$.  A partition is
{\sf $e$-restricted} if $\mu_i-\mu_{i+1}<e$ for $i\ge1$. (This is
compatible with our previous definition of $e$: if $q\ne1$ then~$e$ is
the multiplicative order of $q$ in $R$; otherwise, $e$ is the
characteristic of $R$.)

\begin{Theorem}
[Dipper and James~\cite{DJ:reps}]
Suppose that $R$ is a field. Then $D^\mu\ne0$ if and only if
$\mu$ is $e$-restricted.
\label{A simples}\end{Theorem}

Actually, Dipper and James showed that the simple $\H$--modules are
indexed by $e$-regular partitions (that is, a partition with no~$e$
non-zero parts being equal). Our statement is different from theirs
because our Specht modules are isomorphic to the duals of the
Dipper-James Specht modules~\cite{murphy:basis}. 

Using the $\L$-module structure of the Specht modules it is
straightforward to see that $D^\mu\ne0$ whenever~$\mu$ is
$e$-restricted (recall that $\L=\<L_1,\dots,L_n\>$). The converse is
harder and follows from showing that if~$\mu$ is not $e$-restricted
then $[e]_q!=\prod_{k=1}^e(1+q+\dots+q^{k-1})$ divides the Gram
determinant of the Specht module defined over $\Z[q,q^{-1}]$.  For the
proof see \cite{M:ULect,murphy:basis,DJ:reps}.

For $r>1$ the next result follows easily from Theorem~\ref{A simples}. The
statement is misleading because two separate, but similar, arguments are
needed. For the proof when $q=1$ see Mathas~\cite{aksimples}; for the
case where $Q_s=0$, for all $s$, see Ariki-Mathas~\cite{AM:simples}.

\begin{Corollary}
[Mathas~\cite{aksimples}, Ariki-Mathas~\cite{AM:simples}]
Suppose that $R$ is a field and that either $($i$)~q=1$, or
$($ii$)~Q_1=\dots=Q_r=0$. Let $\mu=\rtuple\mu$ be a
multipartition of~$n$. Then $D^\mu\ne0$ if and only if the following
two conditions are satisfied.
\begin{enumerate}
\item $\mu^{(s)}$ is $e$-restricted for $1\le s\le r$.
\item $\mu^{(s)}=(0)$ whenever $Q_s=Q_t$ for some $t>s$.
\end{enumerate}
\end{Corollary}

In the case $Q_1=\dots=Q_r=0$ the last result simplifies to saying that
$D^\mu\ne0$ if and only if $\mu=\((0),\dots,(0),\mu^{(r)}\)$ for
some $e$-restricted partition $\mu^{(r)}$.

It remains to treat the cases where $q\ne1$ and
$Q_s\ne0$ for all $s$.

Given two nodes $x=(a,b,s)$ and $y=(c,d,t)$ we say that $y$ is
{\sf below} $x$ if either $s<t$, or $s=t$ and $a<c$. Further,
$x\in[\lambda]$ is {\sf removable} if $[\lambda]\setminus\{x\}$ is the
diagram of a multipartition; similarly, $y\notin[\lambda]$ is {\sf addable}
if $[\lambda]\cup\{y\}$ is the diagram of a multipartition. If
$i=\res(x)$ we call $x$ an $i$-node.

An $i$-node $x$ is {\sf normal} if (i) whenever $y$ is a removable
$i$-node below $x$ then there are more removable $i$-nodes between
$x$ and $y$ than there are addable $i$-nodes, and (ii) there are at
least as many removable $i$-nodes below $x$ as addable $i$-nodes
below $x$.  In addition, $x$ is {\sf good} if there are no normal
$i$-nodes above $x$. If $[\mu]=[\lambda]\cup\{x\}$ for some good node
$x$ we write $\lambda\goodarrow\mu$.

\begin{Definition}
A multipartition $\mu$ is {\sf Kleshchev} if either
$\mu=\((0),,\dots,(0)\)$ or $\lambda\goodarrow\mu$ for some Kleshchev
multipartition $\lambda$.
\end{Definition}

The origin of the definition of the Kleshchev multipartitions is that
they are the vertices of the crystal graph of an integrable
$U_v(\widehat{\mathfrak{sl}}_e)$-module. (When $Q_s=q^{a_s}$, for all
$s$, then the Kleshchev multipartitions are the vertices of the crystal
graph of $L_v(\Lambda)$, where $\Lambda=\sum_{s=1}^r\Lambda_{a_s}$. In
general, we take a direct sum of tensor products of crystal graphs in
accordance with Theorem~\ref{Morita}.) There is an edge in the crystal graph
between two Kleshchev multipartitions if~$\lambda\goodarrow\mu$; the
label of the edge is the residue of the node 
in~$[\mu]\setminus[\lambda]$. For more details see
\cite{AM:simples,JMMO}.

When $r=1$ a partition $\mu$ is Kleshchev if and only if $\mu$ is
$e$-restricted; consequently, as it must, the next result agrees with
Theorem~\ref{A simples} when $r=1$.

\begin{Theorem}
[Ariki~\cite{A:class}]
Suppose that $R$ is a field, $q\ne1$, $Q_s\ne0$ for $1\le s\le r$,
and that $\mu$ is a multipartition of $n$. Then $D^\mu\ne0$ if and
only if $\mu$ is a Kleshchev multipartition.
\label{H simples}\end{Theorem}

The first step towards Theorem~\ref{H simples} is to observe that Theorem~\ref{Morita}
allows us to reduce to the crucial case where $q\ne1$ and
$Q_s=q^{a_s}$ for some integers $a_s$ (a different argument is given
in \cite{AM:simples}). Using Theorem~\ref{Ariki}(ii), Ariki~\cite{A:class} is
able to complete the classification of the irreducible $\H$-modules
over $\C$.  To complete the argument, \cite{AM:simples} shows that the
number of simple modules depends only on the integers $a_s$ and the
multiplicative order of $q$ in $R$.

Finally, we remark that by combining these techniques with results of
Ginzburg~\cite{ChrissGinz}, Ariki and the author~\cite{AM:simples}
classified the simple modules of the affine Hecke algebras over an
algebraically closed field of positive characteristic; again, the hard
work is done by Ariki's paper~\cite{Ariki:can}. When $R=\C$ and $q$
is not a root of unity the simple $\Haff_n$-modules were classified by
Zelevinsky~\cite{Zelevinsky:p-adic}; see also \cite{Rogawski}. Kazhdan
and Lusztig~\cite{KL:DelLang} classified the irreducible
$\Haff_n$-modules when $q\in\C^\times$ is a root of unity.

\subsection{The modular branching rules} One of the most significant
results in modular representation theory from the nineties is
Kleshchev's modular branching rule for the symmetric
groups~\cite{Klesh:I,Klesh:II,Klesh:III,Klesh:IV}. Using a streamlined
version of the same techniques Brundan~\cite{Brundan:KM} extended
these results to the Iwahori-Hecke algebra of the symmetric group.
Using completely different methods, Grojnowski~\cite{Groj:control} and
Grojnowski-Vazirani~\cite{GrojVaz} generalized Kleshchev's
modular branching rules to the Ariki-Koike algebras and the affine
Hecke algebra of type~$A$. (Brundan and Kleshchev~\cite{BK:Cliff} have
also applied Grojnowski's arguments to the projective representations
of the symmetric groups.)

Grojnowski was mainly interested in representations of the affine
Hecke algebra $\Haff_n$; however, as remarked in section~(\ref{affine
H}) every irreducible representation of the affine Hecke algebra is an
irreducible representation for a family of Ariki-Koike algebras. He
studies the functors given by induction and restriction (from
$\Haff_n$ to $\Haff_{n\pm1}$), followed by the taking of socles by
analyzing the effect of these functors on the central characters
of~$\Haff_n$. Grojnowski shows that these functors can be described in
terms of the crystal graphs of integral highest weight modules for the
quantum group $U_v(\widehat{\mathfrak{sl}}_e)$;
cf.~Theorem~\ref{Ariki}(i).

\begin{Theorem}
[Grojnowski~\cite{Groj:control},
Grojnowski-Vazirani~\cite{GrojVaz}]\hfil\newline Suppose that $R$ is a
field, $q\ne1$ and $Q_s\ne0$, for $1\le s\le r$.  Then, for each $m$,
there is an $($unknown$)$ permutation $\pi_m$ of the set Kleshchev
multipartitions of $m$ such that if $\mu$ is a Kleshchev
multipartition then
$$\Soc(\Res(D^\mu))\ \ \cong
   \bigoplus_{\pi_{n-1}(\lambda)\goodarrow\pi_n(\mu)} D^\lambda.$$
\label{modular branching}\end{Theorem}

In \cite{GrojVaz} Grojnowski-Vazirani prove that $\Soc(\Res(D^\mu))$
is multiplicity free. In \cite{Groj:control} Grojnowski shows that
there exists a set of irreducible $\H$--modules which are indexed by
the Kleshchev multipartitions and for which the modular branching rule
is given by removing good nodes; unfortunately, Grojnowski does not
give an explicit construction of these simple modules. Conjecturally,
$\pi_m$ is trivial for all $m$. 

Notice that Theorem~\ref{modular branching} implies that there are at
most $e$ direct summands of $\Soc(\Res(D^\mu))$ and that they all
belong to different blocks.

As Grojnowski remarks, the assumption that $q\ne1$ is not essential
and can be removed (at the expense of some additional notation). Du
and Rui~\cite{DuRui:branching} also obtained the modular branching
rule in the special case where $q^dQ_s\ne Q_t$, for $1\le s<t\le r$
and $|d|<n$. In fact, they obtain the stronger result that
$\pi_m=1$, for all $m$, in this case. By Theorem~\ref{Morita} such
Ariki-Koike algebras are Morita equivalent to direct sums of tensor
products of Iwahori-Hecke algebras $\H_q(\Sym_m)$, so in this case the
result can be deduced from Brundan's theorem~\cite{Brundan:KM}
for~$\H_q(\Sym_m)$.

Grojnowski~\cite{Groj:control} is also able to show that the number of
irreducible $\H_n$--modules is equal to the number of Kleshchev
multipartitions of $n$; this gives a more elementary proof of part of
Theorem~\ref{H simples}. Grojnowski is also able to count the number of
irreducible modules of the affine Hecke algebra~$\Haff_n$ over an
arbitrary algebraically closed field. In section~(\ref{block section})
below we discuss the application of Theorem~\ref{modular branching} to
classifying the blocks of~$\H$.

\section{The cyclotomic $q$-Schur algebra}

This chapter introduces the cyclotomic $q$-Schur algebras. These
algebras are defined as endomorphism algebras
$$\Schur(\Lambda)
    =\End_\H\Big(\bigoplus_{\mu\in\Lambda} M^\mu \Big),$$
where $\Lambda$ is a finite set of multicompositions and $M^\mu$ is a
certain $\H$-module. In the special case where $r=1$ the cyclotomic
$q$-Schur algebras are the $q$-Schur algebras of Dipper and
James~\cite{DJ:Schur}; see \cite{Green,Donkin:book,M:ULect,BDK}.  This
was one of the motivations for introducing the cyclotomic $q$-Schur
algebras.

\subsection{Permutation modules.}

We begin by describing the $\H$-modules $M^\mu$. 

A {\sf composition} of $n$ is a sequence
$\sigma=(\sigma_1,\sigma_2,\dots)$ of non-negative integers $\sigma_i$
such that $|\sigma|=\sum_{i\ge1}\sigma_i=n$; we will sometimes write
$\sigma=(\sigma_1,\dots,\sigma_k)$ if~$\sigma_i=0$ for $i>k$. A {\sf
multicomposition} of~$n$ is an ordered $r$-tuple $\mu=\rtuple\mu$ of
compositions with $|\mu^{(1)}|+\dots+|\mu^{(r)}|=n$. 

\begin{Definition}
Suppose that $\mu$ is a multicomposition of $n$. Then
$M^\mu$ is the right ideal $M^\mu=m_\mu\H$ of $\H$.
\end{Definition}

Given a multicomposition $\mu$ let $\vec\mu=\rtuple{\vec\mu}$ be the
multipartition where $\vec\mu^{(s)}$ is the partition obtained by
ordering the parts of the composition $\mu^{(s)}$. It is not hard to
see that $M^{\vec\mu}\cong M^\mu$; indeed, if $d\in\Sym_n$ is a right
coset representative of $\Sym_\mu$ of minimal length such that
$\Sym_{\vec\mu}=d^{-1}\Sym_\mu d$ then $T_dx_{\vec\mu}=x_\mu T_d$;
hence, $T_dm_{\vec\mu}=m_\mu T_d$ and an isomorphism
$M^{\vec\mu}\cong M^\mu$ is given by $h\mapsto T_dh$,
for $h\in M^{\vec\mu}$.

When $r=1$ these modules are nothing more than the induced trivial
representations of the parabolic subalgebra
$\H_q(\Sym_\mu)=\<T_i\mid t_i\in\Sym_\mu\>
                  =\sum_{w\in\Sym_\mu}RT_w$.
More precisely, let~$\mathbf1_\mu$ be the trivial representation
of the subalgebra $\H_q(\Sym_\mu)$; so~$\mathbf1_\mu$ is a free
$R$-module of rank $1$ on which $T_w$ acts as multiplication by
$q^{\len(w)}$ for all $w\in \Sym_\mu$. Then
$$M^\mu\cong
          \mathbf 1_\mu\otimes_{\H_q(\Sym_\mu)}\H_q(\Sym_n).$$
(Note that $\H_q(\Sym_n)$ is free as a right $\H_q(\Sym_\mu)$-module.)

If $r>1$ then, in general, the modules $M^\mu$ are not obviously
induced from subalgebras (except in the case considered by
Shoji~\cite{Sho:Frob}). Even so, the $M^\mu$ behave very much like
permutation modules, so it is not a bad idea to think to them of them
as such.

In order to describe a basis of $M^\mu$ we need to introduce some more
notation. Let $\nr=\set{(i,s)|i\ge1\And 1\le s\le r}$. If
$(i,s), (j,t)$ are elements of $\nr$ write $(i,s)\preceq(j,t)$ if 
either $s<t$, or $s=t$ and $i\le j$.

Let $\mu$ be a multicomposition. Then a $\lambda$-tableau of {\sf type}
$\mu$ is a map $\T\map{[\lambda]}\nr$ such that
$\mu^{(s)}_i=\#\set{x\in[\lambda]|\T(x)=(i,s)},$ for $1\le s\le r$
and all $i\ge1$; we write $\Type(\T)=\mu$. Again, we will think of
a tableau of type $\mu$ as being an $r$-tuple of tableaux. For example,
two tableaux of type $\((3,1),(1^2),(2,1)\)$ are
$$\Tritab(1_1&1_1&1_1\cr2_1|1_2\cr2_2|1_3&1_3\cr2_3) \And
\Tritab(1_1&1_1&1_1&2_1&1_3\cr1_2&2_3|2_2&1_3|2_3),$$ where we write
$i_s$ instead of the ordered pair $(i,s)$.

If $\s$ is a standard $\lambda$-tableau let $\mu(\s)$ be the
tableau of type $\mu$ obtained by replacing each entry $k$ in $\s$ by
$(i,s)$ if $k$ appears in row $i$ of component $s$ of $\t^\mu$ --- as
for multipartitions, we define $\t^\mu$ to be the $\mu$-tableau with
the integers $1,\dots,n$ entered from left to right and then top to
bottom along the rows of the components of $[\mu]$. 

\begin{Definition}
Let $\lambda$ be a multipartition and $\mu$ a
multicomposition. A {\sf semistandard} $\lambda$-tableau 
a $\lambda$-tableau $\T=\rtuple\T$ such that
\begin{enumerate}
\item the entries in each row of $\T$ are
non-decreasing in each component (when ordered by~$\preceq$); and,
\item the entries in each column of $\T$ are strictly increasing
in each component; and,
\item if
$(a,b,c)\in[\lambda]$ and $\T(a,b,c)=(i,s)$ then $s\ge c$.
\end{enumerate}\noindent Let $\SStd(\lambda,\mu)$ be the set of
semistandard $\lambda$-tableaux of type $\mu$ and let 
$\SStd(\lambda,\Lambda)=\bigcup_{\mu\in\Lambda}\SStd(\lambda,\mu)$.
\label{sstd}\end{Definition}

When $r=1$ condition (iii) is redundant and Definition~\ref{sstd}
becomes the familiar definition of semistandard tableaux from the
representation theory of the general linear and symmetric groups.
Write $\comp_\t(k)=s$ if $k$ appears in component~$s$ of~$\t$. For
$r>1$ condition (iii) is unexpected; it has its origin in the
fact~\cite[Prop.~3.23]{DJM:cyc} that if $h\in M^\mu$ and
$h=\sum_{\s,\t}r_{\s\t}m_{\s\t}$ for some $r_{\s\t}\in R$ then
$\comp_\s(k)\le\comp_{\t^\mu}(k)$, for $k=1,\dots,n$. Observe that
$\mu(\s)$ satisfies condition~(iii) if and only if
$\comp_\s(k)\le\comp_{\t^\mu}(k)$ for all $k$.

For example, if $\lambda$ is a multipartition then
$\T^\lambda=\lambda(\t^\lambda)$ is the unique semistandard
$\lambda$-tableau of type $\lambda$. The first of the two tableaux in
the example above is~$\T^\lambda$ for $\lambda=\((3,1),(1^2),(2,1)\)$;
the second tableau there is also semistandard. Finally, let
$\omega=\((0),\dots,(0),(1^n)\)$. Then it is easy to see that the map 
$$\omega\bij{\Std(\lambda)}\SStd(\lambda,\omega);\s\mapsto\omega(\s)
\Number{sstd=std}$$ 
is a bijection between the set of standard $\lambda$-tableaux and
the set of semistandard $\lambda$-tableaux of type $\omega$.
Hereafter, we identity $\Std(\lambda)$ and $\SStd(\lambda,\omega)$
via (\ref{sstd=std}).

\begin{Definition}
Suppose that $\S$ is a semistandard $\lambda$-tableau of 
type~$\mu$ and that $\t$ is a standard $\lambda$-tableau. Define
$$m_{\S\t}=\sum_{\substack{\s\in\Std(\lambda)\\\S=\mu(\s)}}
                   m_{\s\t}.$$
\end{Definition}

The point of all of this notation is the following useful theorem.

\begin{Theorem}
Suppose that $\mu$ is a multicomposition of $n$. Then~$M^\mu$ is free 
as an $R$-module with basis
$$\Set[47]m_{\S\t}|$\S\in\SStd(\lambda,\mu)$ and $\t\in\Std(\lambda)$
                        for some multipartition $\lambda$ of $n$|.$$
\label{M^mu basis}\end{Theorem}

When $r=1$ this result was first proved by Murphy~\cite{murphy:basis};
the general case can be found in Dipper-James-Mathas~\cite{DJM:cyc}.

The proof of this result is straightforward. A small calculation shows
that $m_{\S\t}$ is an element of $M^\mu$. Next, the elements in the
statement of Theorem~\ref{M^mu basis} are linearly independent by Theorem~\ref{Murphy basis}. Finally, if $h\in M^\mu$ then $h$ can be written as a linear
combination of standard basis elements; in turn, these are a linear
combination of the~$m_{\S\t}$.

The importance of this result stems from the following applications.

\begin{Corollary}
Suppose that $\mu$ is a multicomposition of $n$.  Then there
exists a filtration $M^\mu=M_1>M_2>\cdots>M_{k+1}=0$ of~$M^\mu$ such
that
\begin{enumerate}
\item $M_i/M_{i+1}\iso S^{\lambda_i}$ for some multipartition
$\lambda_i$ for $i=1,\dots,k$; and,
\item for each multipartition $\lambda$ the number
of $i$ with $\lambda=\lambda_i$ is equal to the number of
semistandard $\lambda$-tableaux of type $\mu$.
\end{enumerate}
\end{Corollary}

\begin{proof}
[Sketch of proof]
Fixing $\S$ and varying $\t$ in the basis $\{m_{\S\t}\}$
of $M^\mu$ gives a Specht module modulo higher terms.
\end{proof}

For each semistandard $\lambda$-tableau $\S$ of type $\mu$ and each
semistandard $\lambda$-tableau $\T$ of type $\nu$ define
$$m_{\S\T}=\sum_{\substack{\t\in\Std(\lambda)\\\T=\nu(\t)}}m_{\S\t}.$$
By definition, $m_{\S\T}=\sum_{\s\t}m_{\s\t}$ where the sum is over
the standard $\lambda$-tableaux $\s$ and $\t$ such that $\mu(\s)=\S$
and $\nu(\t)=\T$.

\begin{Corollary}
Suppose that $\mu$ and $\nu$ are multicompositions of $n$. Then
$$\Set[48]m_{\S\T}|$\S\in\SStd(\lambda,\mu)$ and $\T\in\SStd(\lambda,\nu)$
                    for some multipartition $\lambda$ of $n$|$$
is a basis of $\H m_\nu\cap m_\mu\H$.
\label{m_ST basis}\end{Corollary}

We are now ready to tackle the cyclotomic $q$-Schur algebras.

\subsection{The semistandard basis theorem} We give a slightly more
general definition for the cyclotomic $q$-Schur algebras than
appeared in \cite{DJM:cyc} in that we allow the set $\Lambda$ to be an
arbitrary finite set of multicompositions. We invite the reader to
check that the arguments from \cite{DJM:cyc} go through without
change.

Extend the dominance ordering $\gedom$ to the set of all
multicompositions; by restriction we consider any set of
multicompositions as a poset.

\begin{Definition}
Suppose that $\Lambda$ is a finite set of multicompositions of $n$.
The {\sf cyclotomic $q$-Schur algebra} is the endomorphism algebra
$$\Schur(\Lambda)=\End_\H\Big(\bigoplus_{\mu\in\Lambda} M^\mu\Big).$$
Let $\Lambda^+
   =\set{\lambda\vdash n|\lambda\gedom\mu\ForSome\mu\in\Lambda}$.
\label{Schur defn}\end{Definition}

We should really write $\Schur(\Lambda)=\Schur_{R,q,\q}(\Lambda)$
since $\Schur(\Lambda)$ depends on $\Lambda$, $R$, $q$ and $\q$.  

Part of the original definition of the cyclotomic $q$-Schur algebras
in~\cite{DJM:cyc} was the requirement that $\Lambda^+\ss\Lambda$.
Following Donkin~\cite{Donkin:book}, we say that $\Lambda$ is {\sf saturated}
if $\Lambda^+\ss\Lambda$.  In analogy with representations of Lie groups,
$\Lambda^+$ should be thought of as the set of dominant weights and
$\Lambda$ the set of weights. Note that $\Lambda^+$ is not
necessarily a subset of $\Lambda$.

Let $\Lambda=\Lambda(d;n)$ be the set of all compositions
$\mu=(\mu_1,\dots,\mu_d)$ of~$n$ of length at most $d$ (so $\mu_i=0$
whenever $i>d$). Then $\Schur_q(d;n)=\Schur\(\Lambda(d;n)\)$ is one of
the $q$-Schur algebras of Dipper and James~\cite{DJ:Schur}.

As an $R$-module we see that 
$$\Schur(\Lambda)=\End_\H\Big(\bigoplus_{\mu\in\Lambda} M^\mu\Big)
         =\bigoplus_{\mu,\nu\in\Lambda} \Hom_\H\Big(M^\nu,M^\mu\Big);$$
so we need to understand the $R$-modules $\Hom_\H(M^\nu,M^\mu)$.

\begin{Proposition}
Suppose that $\mu$ and $\nu$ are multicompositions of~$n$. 
Then an $R$-linear map $\phi\map{M^\nu}M^\mu$ belongs to 
$\Hom_\H(M^\nu,M^\mu)$ if and only if 
$$\phi(m_\nu)
  =\sum_{\substack{\S\in\SStd(\lambda,\mu)\\\T\in\SStd(\lambda,\nu)}}
         r_{\S\T}m_{\S\T}$$
for some $r_{\S\T}\in R$.
\label{hom basis}
\end{Proposition}

\begin{proof}
[Sketch of proof]
If $Q_1,\dots,Q_r$ are invertible elements of $R$ then $\H$ is
a symmetric algebra by Proposition~\ref{symmetric}(ii); therefore, 
$\Hom_\H(M^\nu,M^\mu)$ and $m_\mu\H\cap\H m_\nu$ are canonically
isomorphic $R$-modules (via the map $\phi\mapsto\phi(m_\nu)$), so the
proposition follows by Corollary~\ref{m_ST basis}. 

For the general case, an intricate induction (see
\cite[\Sect5]{DJM:cyc}), which is independent of
Proposition~\ref{symmetric}, shows that the double annihilator
of $m_\mu$, $$\set{x\in\H|xh=0\Whenever m_\mu h=0\ForSome h\in\H},$$
is~$\H m_\mu$. It again follows that
$\Hom_\H(M^\nu,M^\mu)\cong m_\nu\H\cap\H m_\mu$, so we can
complete the proof using the argument of the last paragraph.
\end{proof}

Observe that a map $\phi\in\Hom_\H(M^\nu,M^\mu)$ is completely determined
by~$\phi(m_\nu)$ since $\phi(m_\nu h)=\phi(m_\nu)h$ for all $h\in\H$.
Therefore, we can lift the involution~$*$ of $\H$ to give an involutory
anti-isomorphism of $\Schur(\Lambda)$ by defining
$\phi^*\in\Hom_\H(M^\mu,M^\nu)$ by 
$\phi^*(m_\mu h)=\(\phi(m_\nu)\)^*h$ for all $h\in\H$.

\begin{Definition}
Suppose that $\lambda\in\Lambda^+$ is a multipartition and that 
$\mu, \nu\in\Lambda$ are multicompositions. For each pair of standard
$\lambda$-tableaux $\S\in\SStd(\lambda,\mu)$
and $\T\in\SStd(\lambda,\nu)$ let $\phiST$ be the
$R$-linear endomorphism of $\bigoplus_{\mu\in\Lambda}M^\mu$
determined by
$$\phiST(m_\alpha h)=\delta_{\alpha\nu}m_{\S\T}h,$$
for all $\alpha\in\Lambda$ and $h\in\H$ (here
$\delta_{\alpha\nu}$ is the Kronecker delta).
\end{Definition}

By Proposition~\ref{hom basis} $\phiST$ is an element of
$\Schur(\Lambda)$. 

Let $\mathcal S^\lambda$ be the $R$-submodule of $\Schur(\Lambda)$
spanned by the $\phi_{\U\V}$, for some $\U,\V\in\SStd(\rho,\Lambda)$
where $\rho\in\Lambda^+$ and~$\rho\gdom\lambda$. From the definitions,
$\mathcal S^\lambda$ consists of those elements of
$\Schur(\Lambda)$ whose image is contained in $\mathcal H^\lambda$. 

We can now state the semistandard basis theorem for the cyclotomic
Schur algebras.

\begin{Theorem}[%
\protect{Dipper-James-Mathas~\cite[Theorem~6.6]{DJM:cyc}}]
\label{SSbasis}\hfil\newline
Let $\Lambda$ be a finite set of multicompositions. Then
the cyclotomic $q$-Schur algebra $\Schur(\Lambda)$ is free as an 
$R$-module with basis
$$\set{\phiST|\ForSome \S,\T\in\SStd(\lambda,\Lambda)\And
            \lambda\in\Lambda^+}.$$
Moreover, this basis is a cellular basis of $\Schur(\Lambda)$; more
precisely, if $\S$ and~$\T$ are semistandard $\lambda$-tableaux, for
some $\lambda\in\Lambda^+$, then
\begin{enumerate}
\item $\phi_{\S\T}^*=\phi_{\T\S}$; and,
\item for all $\phi\in\Schur(\Lambda)$ there exist scalars 
$r_\V=r_{\T\V}(\phi)\in R$, which do not depend on~$\S$, such that
$$\phiST\phi\equiv\sum_{\V\in\SStd(\lambda,\Lambda)}
             r_\V\phi_{\S\V}\pmod{\mathcal S^\lambda}.$$
\end{enumerate}\end{Theorem}

\begin{proof}
[Sketch of proof]
Proposition~\ref{hom basis} implies that these elements give a basis of
$\Schur(\Lambda)$. Using Theorem~\ref{Murphy basis} it is not hard to
see that the semistandard basis is cellular. 
\end{proof}

In particular, notice that $\Schur(\Lambda)$ is always free as an
$R$-module and that its rank is independent of $R$, $q$ and $\q$.
The semistandard basis of $\Schur(\Lambda)$ really comes from
Theorem~\ref{M^mu basis} and the basis element $\phiST$ really comes
from a Specht filtration of $M^\mu$.

It is worthwhile explaining how the multiplication in
$\Schur(\Lambda)$ is determined. Suppose that $\S,\T,\U$ and $\V$ are
semistandard tableaux and suppose that $\nu=\Type(\V)$ and
$\mu=\Type(\U)$. Then $m_{\U\V}=m_\mu h_{\U\V}$, for some
$h_{\U\V}\in\H$, and
$$\phi_{\S\T}\phi_{\U\V}=\sum_{\sA,\sB}r_{\sA\sB}\phi_{\sA\sB},$$
where the scalars $r_{\sA\sB}\in R$ are determined by
$m_{\S\T}h_{\U\V}=\sum r_{\sA\sB}m_{\sA\sB}$; this makes sense by
Proposition~\ref{m_ST basis} and is proved by evaluating the functions on
both sides at $m_\nu$. Note, in particular, that $r_{\sA\sB}=0$ unless
$\Type(\U)=\Type(\T)$, $\Type(\sA)=\Type(\S)$ and
$\Type(\sB)=\Type(\V)$. In Theorem~\ref{SSbasis}(ii), $r_\V=r_{\S\V}$.

With some work it is possible to show that when $r=1$ this basis
agrees with Richard Green's codeterminant basis of the $q$-Schur
algebra~\cite{RGreen:codet}; see also \cite{Green:codet,Wood:codet}.
When $r=2$ and $\H$ is symmetric Theorem~\ref{SSbasis} is equivalent to a
theorem of Du and Scott~\cite{DuScott:BSchur}.

\subsection{Weyl modules for cyclotomic $q$-Schur algebras} 
By the semistandard basis theorem $\Schur(\Lambda)$ is a cellular
algebra. Therefore, exactly as in Definition~\ref{Specht modules} we
can write down a collection of {\it cell modules} for
$\Schur(\Lambda)$ and, up to isomorphism, every irreducible
$\Schur(\Lambda)$-module is a quotient of one of these modules.

\begin{Definition}
Suppose that $\lambda\in\Lambda^+$ is a multipartition. The
{\sf Weyl module} $W^\lambda$ is the free $R$-module with basis
$\set{\phi_\T|\T\in\SStd(\lambda,\Lambda)}$ on which
$\phi\in\Schur(\Lambda)$ acts via
$$\phi_\T\phi=\sum_{\V\in\SStd(\lambda,\Lambda)}r_\V\phi_\V,$$
where the scalars $r_\V\in R$ are as in Theorem~\ref{SSbasis}(ii).
\end{Definition}

It follows from Theorem~\ref{SSbasis} that $W^\lambda$ is a right
$\Schur(\Lambda)$-module. As with the Specht modules we define a
bilinear form on $W^\lambda$ by requiring that 
$\<\phi_\S,\phi_\T\>\phi_{\U\V}
    \equiv\phi_{\U\S}\phi_{\T\V}\pmod{\mathcal S^\lambda}$
for semistandard tableaux $\S,\T,\U,\V\in\SStd(\lambda,\Lambda)$. 
Then the radical of this form, $\Rad W^\lambda$, is a submodule of
$W^\lambda$ and we define $L^\lambda=W^\lambda/\Rad W^\lambda$.

Exactly as in Theorem~\ref{main H theorem}, the theory of cellular algebras
now gives us the following.

\begin{Theorem}
\label{S modules} Suppose that $R$ is a field.
\begin{enumerate}
\item For each $\lambda\in\Lambda^+$,
$L^\lambda$ is either zero or an absolutely
irreducible $\Schur(\Lambda)$-module.
\item $\set{L^\lambda|\lambda\in\Lambda^+\And L^\lambda\ne0}$ is a 
complete set of pairwise non-isomorphic irreducible 
$\Schur(\Lambda)$-modules.
\item $\Schur(\Lambda)$ is semisimple if and only if
$L^\lambda=W^\lambda$ for all $\lambda\in\Lambda^+$.
\item Suppose that $\mu,\lambda\in\Lambda^+$ and $L^\lambda\ne0$. Then 
$[W^\mu{:}L^\lambda]\ne0$ only if $\mu\gedom\lambda$; moreover,
$[W^\lambda{:}L^\lambda]=1$.
\end{enumerate}\end{Theorem}

At this level of generality, determining exactly when $L^\lambda$ is 
non-zero is a difficult task. To see this notice that if
$\Lambda=\{\omega\}$ then $\Schur(\Lambda)=\End_\H(\H)\cong\H$ and
$\Lambda^+$ is the set of all partitions of~$n$; so Theorem~\ref{H
simples} is a special case of Theorem~\ref{S modules}. When the poset
$\Lambda$ is saturated (that is, $\Lambda^+\ss\Lambda$) we can say
much more.

Assume now that $\Lambda^+\ss\Lambda$ and let $\lambda$ be a
multipartition of $n$. Then $M^\lambda$ is a summand
of $\oplus_{\mu\in\Lambda}M^\mu$ and so the identity map
$\phi_\lambda\map{M^\lambda}M^\lambda$ is an element of
$\Schur(\Lambda)$. Indeed, looking at the definitions,
$\phi_\lambda=\phi_{\T^\lambda\T^\lambda}$, where
$\T^\lambda=\lambda(\t^\lambda)$ is the unique semistandard
$\lambda$-tableau of type $\lambda$. It follows that the Weyl module
$W^\lambda$ is isomorphic to the submodule of
$\Schur(\Lambda)/\mathcal S^\lambda$ generated by
$\phi_\lambda+\mathcal S^\lambda$, the isomorphism being given by
$$\phi_\T\mapsto\phi_{\T^\lambda\T}+\mathcal S^\lambda
     =(\phi_\lambda+\mathcal S^\lambda)\phi_{\T^\lambda\T},$$ 
for all $\T\in\SStd(\lambda,\Lambda)$.

\begin{Theorem}
\label{quasi}
Suppose that $R$ is a field and that $\Lambda^+\ss\Lambda$.
\begin{enumerate}

\item $L^\lambda$ is a non-zero absolutely irreducible
$\Schur(\Lambda)$-module for all $\lambda\in\Lambda^+$.
\item $\Schur(\Lambda)$ is a quasi-hereditary algebra.
\end{enumerate}\end{Theorem}

\begin{proof}
[Sketch of proof] To prove (i) observe that
$\phi_\lambda=\phi_{\T^\lambda\T^\lambda}$ is an element of
$\Schur(\Lambda)$ because $\lambda\in\Lambda$. Therefore,
$\phi_{\T^\lambda}\in W^\lambda$ and so
$$\<\phi_{\T^\lambda},\phi_{\T^\lambda}\>\phi_{\T^\lambda\T^\lambda} 
   \equiv \phi_{\T^\lambda\T^\lambda}\phi_{\T^\lambda\T^\lambda}
   =\phi_{\T^\lambda\T^\lambda}\pmod{\mathcal S^\lambda};$$ 
hence, $\<\phi_{\T^\lambda},\phi_{\T^\lambda}\>=1$ and
$\phi_{\T^\lambda}\notin\Rad W^\lambda$; so $L^\lambda\ne0$. Part (ii)
follows from~(i) and the structure of cellular algebras.
\end{proof}

Parshall and Wang~\cite{PW} were the first to show that the $q$-Schur
algebras are quasi-hereditary. More generally, the argument above
shows that the $q$-Schur algebras and the cyclotomic Schur algebras
are integrally quasi-hereditary in the sense of~\cite{DPS:qWeylRecip}.

As the example $\Lambda=\{\omega\}$ indicates, when $\Lambda$ is not
saturated the classification of the simple $\Schur(\Lambda)$-modules is
non-trivial. Nor are there obvious necessary and sufficient conditions
for when $\Schur(\Lambda)$ is quasi-hereditary. The answers to these
questions will depend on $\Lambda$, $R$ and the parameters
$q,Q_1,\dots,Q_r$.

The final result of this section is the analogue of
Theorem~\ref{Morita} for the cyclotomic Schur algebras. In
\cite{DM:Morita} a general version of the result below is proved for
an arbitrary finite set of (saturated) multicompositions; we state
only a special case in order to avoid introducing additional notation.

Let $\Lambda_{r,n}$ be the set of all multicompositions of $n$ of
length at most~$n$ and let $\Lambda_{r,n}^+\ss\Lambda_{r,n}$ be the
set of multipartitions of $n$. We write
$\Schur(\Lambda_{r,n})=\Schur_{q,\q}(\Lambda_{r,n})$ to emphasize the
choice of parameters. 

\begin{Theorem}
[Dipper-Mathas~\cite{DM:Morita}]\label{SMorita}
Suppose that $R$ is an integral domain and let 
$\q=\q_1\coprod\dots\coprod\q_\kappa$ be a partition
of $\q$ and suppose that the polynomial $P_n(q,\q_1,\dots,\q_\kappa)$
is invertible in $R$. Let $r_\alpha=|\q_\alpha|$, for
$1\le\alpha\le\kappa$. Then $\Schur(\Lambda_{r,n})$ is Morita
equivalent to the $R$-algebra 
$$\bigoplus_{\substack{n_1,\dots,n_\kappa\ge0\\n_1+\dots+n_\kappa=n}}
\Schur_{q,\q_1}(\Lambda_{r_1,n_1})\otimes\dots
          \otimes\Schur_{q,\q_\kappa}(\Lambda_{r_\kappa,n_\kappa}).$$
\end{Theorem}

This result is deduced from Theorem~\ref{Morita} using the theory of Young
modules for Ariki-Koike algebras~\cite{M:tilting}. 

\section{The representation theory of cyclotomic $q$-Schur algebras}
This chapter gives a summary of the main results in the
representation theory of the cyclotomic $q$-Schur algebras. All of
these results are generalizations of theorems for the $q$-Schur
algebras. Throughout we assume that $\Lambda$ is saturated; that is,
$\Lambda^+\ss\Lambda$.

\subsection{A Schur functor and double centralizer property}

Throughout this section we assume that $\omega\in\Lambda$; because of
this $\phi_\omega=\phi_{\T^\omega\T^\omega}$ is an element of
$\Schur(\Lambda)$. Now, $\phi_\omega$ is the identity map on $\H$; in
particular, it is an idempotent. Moreover, it is easy to see that
$\H\cong\phi_\omega\Schur(\Lambda)\phi_\omega$. Hence, by general
nonsense (see, for example \cite{Green,BDK}), $\phi_\omega$ gives rise
to a functor $\Phi_\omega$ from the category of
$\Schur(\Lambda)$-modules to the category of $\H$-modules;
explicitly, if $M$ is a right $\Schur(\Lambda)$-module then
$\Sfun(M)=M\phi_\omega$ is a right $\H$-module.

Notice that the condition $\omega\in\Lambda$ is the analogue for the
cyclotomic Schur algebras of the familiar requirement that $d\ge n$
for the $q$-Schur algebra $\Schur_q(d;n)$.

{\samepage
\begin{Theorem}
[The cyclotomic Schur functor~\cite{JM:cyc-Schaper}] 
Suppose that $R$ is a field and that $\omega\in\Lambda^+\ss\Lambda$.
Let $\lambda\in\Lambda^+$.  Then, as right $\H$-modules,
\begin{enumerate}
\item $\Sfun(W^\lambda)\iso S^\lambda$;
\item $\Sfun(L^\lambda)\iso D^\lambda$.
\end{enumerate}\noindent 
Furthermore, if $D^\mu\ne0$ then 
$[W^\lambda{:}L^\mu]=[S^\lambda{:}D^\mu]$.
\label{Schur functor}\end{Theorem}}

\begin{proof}
[Sketch of proof]
This can be proved either by general arguments as in \cite{Green}.
Alternatively, from the definitions and the semistandard basis theorem
it is clear that $\Sfun(W^\lambda)=S^\lambda$ (if
$\T\in\SStd(\lambda,\mu)$ then 
$\phi_\T\phi_\omega=\delta_{\mu\omega}\phi_\T$). Next observe that if
$\s$ and $\t$ are standard tableaux then the definition of the inner
product on $W^\lambda$ is that 
\begin{align*}
\<\phi_\s,\phi_\t\>\phi_\lambda
  &\equiv\phi_{\t^\lambda\s}\phi_{\t\t^\lambda}\pmod{\mathcal S^\lambda}.
\intertext{Evaluating the functions on both sides at $m_\lambda$ 
we find that}
\<\phi_\s,\phi_\t\>m_\lambda
   &\equiv m_{\t^\lambda\s}m_{\t\t^\lambda}
    \equiv \<m_\s,m_\t\>m_\lambda\pmod{\mathcal H^\lambda}.
\end{align*}
Hence, $\<\phi_\s,\phi_\t\>=\<m_\s,m_\t\>$ and the remaining claims 
follow.
\end{proof}

An important consequence of Theorem~\ref{Schur functor} is that the
decomposition matrix of $\H$ is a submatrix of the decomposition
matrix of $\Schur(\Lambda)$.

\begin{Corollary}
Suppose that $R$ is a field and that
$\omega\in\Lambda^+\ss\Lambda$. Then the decomposition matrix of $\H$
is the submatrix of the decomposition matrix of $\Schur(\Lambda)$ obtained
by deleting those columns indexed by the multipartitions $\mu$ such
that $D^\mu=0$.
\label{submatrix}\end{Corollary}

Observe that
$\bigoplus_{\lambda\in\Lambda}M^\lambda$ is an
$(\Schur(\Lambda),\H)$-bimodule. In fact, each algebra is the full
centralizer algebra for the other and we have a cyclotomic analogue of
Schur-Weyl duality.

\begin{Theorem}
[Double centralizer property]\quad Suppose that
$\omega\in\Lambda$ and that $\Lambda^+\ss\Lambda$.  Then
$$\Schur(\Lambda)\iso\End_\H\(\bigoplus_{\lambda\in\Lambda}M^\lambda\)
\quad\And\quad
\H\iso\End_{\Schur(\Lambda)}\(\bigoplus_{\lambda\in\Lambda}M^\lambda\).
$$
\label{double centralizer}
\end{Theorem}

\begin{proof}
[Sketch of proof]
The first isomorphism is just the definition of 
$\Schur(\Lambda)$ so there is nothing to prove here. For the second
isomorphism for each $\lambda\in\Lambda$ let $\phi_\lambda$ be the
identity map on~$M^\lambda$ and let $\mathscr
M^\lambda=\phi_\lambda\Schur(\Lambda)$. (So~$\mathscr M^\lambda$ is an
$\Schur(\Lambda)$--module and $M^\lambda$ is an $\H$--module.) Then 
there an isomorphism of $\H$-modules
$$\bigoplus_{\lambda\in\Lambda}M^\lambda
    \cong\bigoplus_{\lambda\in\Lambda}\Sfun(\mathscr M^\lambda)
    =\bigoplus_{\lambda\in\Lambda}\phi_\lambda\Schur(\Lambda)\phi_\omega.$$
By definition $\sum_\lambda\phi_\lambda$ is the identity of
$\Schur(\Lambda)$, so
$\Schur(\Lambda)=\bigoplus_\lambda\phi_\lambda\Schur(\Lambda)$ and
$\bigoplus_{\lambda\in\Lambda}M^\lambda\cong\Schur(\Lambda)\phio$ as 
a left $\Schur(\Lambda)$-module. Therefore,
$$\End_{\Schur(\Lambda)}\(\bigoplus_{\lambda\in\Lambda}M^\lambda\)
       \iso\End_{\Schur(\Lambda)}\(\Schur(\Lambda)\phio\)
       \iso\phio\Schur(\Lambda)\phio.$$
As $\phio\Schur(\Lambda)\phio\cong\H$, this completes the proof.
\end{proof}

\subsection{The blocks of the cyclotomic Schur algebras}
\label{block section}
The centre of the affine Hecke algebra $\Haff_n$ is given by the
following well-known result of Bernstein.

\begin{Theorem}
[Bernstein]\label{Bernstein}
Suppose that $R$ is an algebraically closed field. Then the centre of
$\Haff_n$ is equal to $R[X_1^\pm,\dots,X_n^\pm]^{\Sym_n}$, the
$R$-algebra of symmetric Laurent polynomials in~$X_1,\dots,X_n$.
\end{Theorem}

This is quite straightforward to prove given the Bernstein
presentation of $\Haff_n$. 

Now $X_k$ maps to $L_k$ under the natural surjection $\Haff_n\To\H_n$ so
this implies that any symmetric polynomial in $L_1,\dots,L_n$ belongs
to the centre of the Ariki-Koike algebra $\H_n$. As we remarked
earlier, in the semisimple case the centre of $\H_n$ is always the
algebra of symmetric polynomials in $L_1,\ldots,L_n$; however, when
$\H_n$ is not semisimple the centre of~$\H_n$ can be larger than this.
Because of this Theorem~\ref{blocks} below is a little surprising. First, some
notation.

Given a multipartition $\lambda$ let
$\res(\lambda)=\set{\res_{\t^\lambda}(k)|1\le k\le n}$, which we
consider as a {\it multiset}. By the remarks above and
Proposition~\ref{L action} if two simple $\H_n$-modules $D^\lambda$ and
$D^\mu$ are in the same block then $\res(\lambda)=\res(\mu)$ as
multisets.

We also note that because $\H_n$ is a cellular algebra all of the
composition factors of $S^\lambda$ belong to the same block; hence,
$D^\lambda$ and $D^\mu$ are in the same block if and only  if
$S^\lambda$ and~$S^\mu$ are in the same block. The same remark applies
to the simple modules and the Weyl modules of the cyclotomic
$q$-Schur algebras.

\begin{Theorem}
\label{blocks} 
Suppose that $R$ is an algebraically closed field and that $\lambda$
and $\mu$ are multipartitions of $n$. Then the following are
equivalent.
\begin{enumerate}
\item $\res(\lambda)=\res(\mu)$ as multisets.
\item $S^\lambda$ and $S^\mu$ are in the same block as $\H_n$-modules.
\item $S^\lambda$ and $S^\mu$ are in the same block as 
$\Haff_n$-modules.
\item $W^\lambda$ and $W^\mu$ are in the same block as
$\Schur(\Lambda_{r,n})$-modules.
\end{enumerate}\end{Theorem}

\begin{proof}
[Sketch of proof]
As noted above, the implication (ii)$\Rightarrow$(i) follows from
Proposition~\ref{L action}; this is was first proved by Graham and
Lehrer~\cite{GL} who also conjectured that the converse was true. That
(i) and (iii) are equivalent follows from Theorem~\ref{Bernstein}. 

The hard part is proving that (iii) implies (ii); this was done by
Grojnowski~\cite{Groj:AKblocks} using his modular branching rule.  The
key point is that if $\lambda$ and $\mu$ are distinct multipartitions
with $D^\lambda\ne0$, $D^\mu\ne0$ and $\res(\lambda)=\res(\mu)$ then 
$\Hom_{\Haff_{n-1}}(\Res D^\lambda,\Res D^\mu)=0$ by 
Theorem~\ref{modular branching}; here $\Res$ is the functor
$\Res\map{\Haff_n\mod}\Haff_{n-1}\mod$. Grojnowski shows that this
implies that whenever $0\To D^\lambda\To X\To D^\mu\To 0$ is an exact
sequence of $\Haff_n$-modules then it is still exact when considered
as a sequence of $\H_n$-modules (for any $\H_n$-module $X$). This
implies (ii).

Finally, by the double centralizer property (Theorem~\ref{double
centralizer}), $\Schur(\Lambda_{r,n})$ and $\H_n$ have the same number
of blocks (see \cite[Cor.~5.38]{M:ULect}), so it follows that (ii) and
(iv) are equivalent.
\end{proof}


Theorem~\ref{blocks} does not classify the blocks of an arbitrary
cyclotomic Schur algebra; rather it classifies the blocks of
$\Schur(\Lambda)$ for any $\Lambda$ with $\Lambda^+_{r,n}\ss\Lambda$
(by standard arguments, all of these algebras are Morita equivalent).
When $r=1$ the blocks for the $q$-Schur algebras $\Schur_q(d;n)$ have
been classified by Cox~\cite{Cox:blocks} for all $d$, $n$ and~$q$.
The general case is open when $r>1$.

\subsection{The Jantzen sum formula} Throughout this section assume
that $R$ is a field and that $\Lambda$ is saturated. Let~$t$ be an
indeterminate over $R$ and let $\p$ be the maximal ideal of
$R[t,t^{-1}]$ generated by~$t-1$. The localization
$\O=R[t,t^{-1}]_{\p}$ of $R[t,t^{-1}]$ at $\p$ is a discrete valuation
ring and $R\cong\O/\p$.  Let~$\nu_\p$ be the $\p$-adic valuation
on~$\O$. 

Let $\H_\O$ be the Hecke algebra over $\O$ with parameters $qt$ and
$U_s=Q_st^{ns}$ if $Q_s\ne0$ and $U_s=(t^{ns}-1)$ if $Q_s=0$. Then 
$\H_{R(t)}=\H_\O\otimes R(t)$ is semisimple by Corollary~\ref{H ss} and
$\H_R=\H_{q,\q}(W_{r,n})\cong\H_\O\otimes_\O R$ is the reduction of
$\H_\O$ modulo $\p$. Let $\Schur_\O(\Lambda)$, $\Schur_{R(t)}(\Lambda)$
and $\Schur_R(\Lambda)$ be the corresponding cyclotomic $q$-Schur
algebras.

Define the {\sf $\O$-residue} of a node $x=(i,j,s)$ to be
$\res_\O(x)=(qt)^{j-i}U_s$, an element of~$\O$. The connection with
our previous definition of residue is that 
$\res(x)=\res_\O(x)\otimes_\O 1_R$.

Let $\lambda$ be a multipartition and for each node
$x=(i,j,s)\in[\lambda]$ let $r_x\ss[\lambda]$ be the corresponding rim
hook (so $r_x$ is a rim  hook in $[\lambda^{(s)}]$); then
$[\lambda]\setminus r_x$ is the diagram of a multipartition. Let
$\llen(r_x)$ be the leg length of $r_x$ and define
$\res_\O(r_x)=\res_\O(f_x)$ where~$f_x$ is the foot node of~$r_x$.
These definitions can be found in \cite{JK,M:ULect}.

Suppose that $\lambda$ and $\mu$ are multipartitions of $n$. If
$\lambda\not\gdom\mu$ let $g_{\lambda\mu}=1$; otherwise set
$$g_{\lambda\mu}=\prod_{x\in[\lambda]} 
      \prod_{\SR{y\in[\mu]}{[\mu]\setminus r_y=[\lambda]\setminus r_x}}
             \(\res_\O(r_x)-\res_\O(r_y)\)^{\varepsilon_{xy}},$$
where $\varepsilon_{xy}=(-1)^{\llen(r_x)+\llen(r_y)}$. The scalars
$g_{\lambda\mu}\in\O$ have a combinatorial interpretation in terms of
moving rim hooks in the diagram of multipartitions; see
\cite[Example~3.39]{JM:cyc-Schaper}.

Finally, let $W^\lambda_\O$ and $W^\lambda_R$ be the Weyl modules
for $\Schur_\O(\Lambda)$ and $\Schur_R(\Lambda)$ respectively; note that
$W^\lambda_R\cong W^\lambda_\O\otimes_\O R$ as $R$-modules. For each
$i\ge0$ define 
$$W^\lambda_\O(i)
  =\set{x\in W^\lambda_\O|\<x,y\>\in\p^i\ForAll y\in W^\lambda_\O}$$
and set
$W^\lambda_R(i)
      =\(W^\lambda_\O(i)+\p W^\lambda_\O\)/\p W^\lambda_\O$. 
The {\sf Jantzen filtration} of $W^\lambda_R$ is 
$$W^\lambda_R=W^\lambda_R(0)\ge W^\lambda_R(1)
              \ge W^\lambda_R(2)\ge\cdots.$$
In particular, note that $\Rad W^\lambda_R=W^\lambda_R(1)$;
consequently, $W^\lambda_R(1)$ is a proper submodule 
of~$W^\lambda_R(0)$ and 
$W^\lambda_R(0)/W^\lambda_R(1)\cong L^\lambda_R$. Note also that
$W^\lambda_R(k)=0$ for $k\gg0$.

Actually, what we have just given is a special case of the definition
of a Jantzen filtration. More generally, the same construction gives a
Jantzen filtration for any suitable modular system $(K,\O,\p)$ (with
parameters). The point of this remark is that the Jantzen filtration
of $W^\lambda_R$ depends upon a non-canonical choice of modular system.

We can now state the analogue of Jantzen's sum formula for
$\Schur_R(\Lambda)$.

\begin{Theorem}
[\protect{James-Mathas~\cite[Theorem~4.6]{JM:cyc-Schaper}}]
Let $\lambda$ be a multipartition of $n$. Then
$$\sum_{i>0}[W^\lambda_R(i)]
   =\sum_{\mu:\lambda\gdom\mu}\nu_\p(g_{\lambda\mu})[W^\mu_R].$$
in the Grothendieck group $K_0(\Schur_R(\Lambda)\text{\mod})$ of 
$\Schur_R(\Lambda)$.
\label{sum formula}\end{Theorem}

When $r=1$ this result describes the Jantzen filtration of the Weyl
modules of the $q$-Schur algebra. The Weyl modules of the $q$-Schur
algebra coincide with the Weyl modules of quantum $\gl_d$; therefore,
when $r=1$ Theorem~\ref{sum formula} is a special case of a result of
Andersen, Polo and Wen~\cite{APW} who proved the analogue of the
Jantzen sum formula for the quantum groups of finite type as a
consequence of Kemp's vanishing theorem. For a combinatorial proof
which takes place inside the $q$-Schur algebra
see~\cite{JM:Schaper,M:ULect}. When $r>1$ there is no geometry to work
with. The argument of \cite{JM:cyc-Schaper} generalizes that of
\cite{JM:Schaper}.

The idea behind the proof of Theorem~\ref{sum formula} is quite
simple.  First, for each $\mu$ compute the determinant of the Gram
matrix $G^\lambda_\mu =(\<\phi_\S,\phi_\T\>)$,
${\S,\T\in\SStd(\lambda,\mu)}$, of the $\mu$-weight space
$W^\lambda_\O\phi_\mu$ of $W^\lambda_\O$. It turns out that 
$\det G^\lambda_\mu=g_{\lambda\mu}$. Now, the inner product 
$\<\ ,\ \>$ on $W^\lambda_\O$ is non-degenerate; so Jantzen's elementary,
yet fundamental, lemma says that 
$$\sum_{i>0}\dim_R W^\lambda_R(i)\phi_\mu =\nu_\p(\det G^\lambda_\mu).$$ 
This is enough to deduce the result because, by Theorem~\ref{S
modules}(iv), any $\Schur(\Lambda)$-module is uniquely determined by
the dimensions of its weight spaces since $\dim L^\nu\phi_\nu=1$ for
all $\nu\in\Lambda^+$.  

Of course, computing $\det G^\lambda_\mu$ is not so easy. This is
accomplished using an orthogonal basis of $W^\lambda_R$ when
$P_\H(q,\q)\ne0$. With this basis the Gram determinant is easier to
calculate because almost all inner products are zero (we are really
computing inner products in $W^\lambda_{R(t)}$). The orthogonal basis
is constructed using a family of operators which act in a triangular
fashion on the semistandard basis of the Weyl modules; intuitively,
these operators belong to something like a Cartan subalgebra
of~$\Schur(\Lambda)$ ---~in fact, they are `lifts' of the elements
$L_k$ to $\Schur(\Lambda)$.

The definition of the Jantzen filtration only requires a finitely
generated $\O$-module which possesses a non-degenerate bilinear form.
The same construction gives a Jantzen filtration 
$S^\lambda=S^\lambda_R(0)\ge S^\lambda_R(1)\ge\cdots$
for each Specht module; equivalently, by the proof of 
Theorem~\ref{Schur functor} we can set
$S^\lambda_R(i)=\Sfun(W^\lambda_R(i))$. Applying the Schur functor
to Theorem~\ref{sum formula} yields the following.

\begin{Corollary}[James-Mathas~\cite{JM:cyc-Schaper}] 
Let $\lambda$ be a multipartition of~$n$. Then
$$\sum_{i>0}[S^\lambda_R(i)]
   =\sum_{\mu:\lambda\gdom\mu}\nu_\p(g_{\lambda\mu})[S^\mu_R].$$
in the Grothendieck group $K_0(\H_R\text{\mod})$ of $\H_R$.
\end{Corollary}

For the symmetric groups (that is, $r=1$ and $q=1$) this is a result
of long standing known as Schaper's Theorem~\cite{Schaper}.  Schaper's
argument is a translation of the Jantzen sum formula for the Weyl
modules of the general linear group~\cite{Jant:kon} (phrased in terms
of the dot action of the  symmetric group upon the weight lattice
of~$\GL_n$), into the combinatorial language of the symmetric group.
It is worth remarking that the corresponding result for the Weyl
groups of type $B$ (i.e.,~$r=2$ and $q=1$), was obtained only
relatively recently \cite{JM:cyc-Schaper}.

The main application of the cyclotomic sum formula has been a
classification of the irreducible Weyl modules and the irreducible
Specht modules with $S^\lambda=D^\lambda$; see \cite{JM:cyc-Schaper}.
When $r=1$ the sum formula was used to complete the classification of
the blocks of the $q$-Schur algebras and the Iwahori-Hecke algebras
of type $A$ and also to classify the ordinary irreducible
$\GL_n(q)$-modules which remain irreducible when reduced mod~$p$ when
$p\nmid q$; see~\cite{JM:Schaper}.  Ariki and the
author~\cite{AM:reptype} have also used the Jantzen sum formula to
classify the representation type of the Iwahori-Hecke algebras of
type~$B$.

\subsection{Connections with quantum groups}\label{SS}
For this section only we renormalize the basis of the Ariki-Koike
algebras so as to be consistent with the notation in
\cite{Ariki:Schur,SakSho}.  We assume that $q$ has a square root in
$R$ and let $q=v^2$. As every field is a splitting field for $\H$ and
$\Schur(\Lambda)$ we are free to extend $R$ so that it contains a
square root of~$q$ if necessary.

Let $\tilde T_i=v^{-1}T_i$ for $1\le i<n$. Then 
$T_0,\tilde T_1,\dots,\tilde T_{n-1}$ still generate $\H$ and they are
subject to the same relations as before except that the quadratic
relation for the $T_i$ becomes $(\tilde T_i-v)(\tilde T_i+v^{-1})=0$,
for $1\le i<n$. Observe that $L_k=\tilde T_{k-1}\dots\tilde
T_1T_0\tilde T_1\dots\tilde T_{k-1}$ for $k=1,\dots,n$.

Fix an integer $d\ge1$ and let $U_v(\gl_d)$ be the quantized
enveloping algebra of $\gl_d$. Thus, $U_v(\gl_d)$ is an associative
$\Q(v)$-algebra which is generated by elements $E_i, F_i, K_j^\pm$,
where $1\le i<n$ and $1\le j\le n$, which are subject to the quantum
Serre relations.

Let $V$ be a $d$ dimensional $\Q(v)$-vector space with basis
$\{e_1,\dots,e_d\}$. Then $V$ is naturally a $U_v(\gl_d)$-module, 
where the action of $U_v(\gl_d)$ on~$V$ is determined by
$$E_ie_a=\delta_{a,i+1}e_{a-1},\quad
F_ie_a=\delta_{a,i}e_{a+1}, \quad\And\quad
K_j e_a=v^{\delta_{j,a}}e_a,$$
for $1\le i<n$, $1\le j\le n$ and $1\le a\le n$.  Now, $U_v(\gl_d)$
is a Hopf algebra with coproduct $\Delta$ given by
$\Delta(K_j)=K_j\otimes K_j$,
$$\Delta(E_i)=E_i\otimes K_i+1\otimes E_i
\quad\And\quad
  \Delta(F_i)=F_i\otimes 1+K_i^{-1}\otimes F_i,$$
for $1\le i<n$ and $1\le j\le n$. Therefore, $V^{\otimes n}$ is a
$U_v(\gl_d)$-module; let $\rho_n\map{U_v(\gl_d)}\End(V^{\otimes n})$
be the corresponding representation.

Let $I(d;n)=\set{(a_1,\dots,a_n)|1\le a_1,\dots,a_n\le d}$. If
$\a\in I(d;n)$ let 
$e_{\a}=e_{a_1}\otimes\dots\otimes e_{a_n}$. Then 
$\set{e_{\a}|\a\in I(d;n)}$ is a basis of $V^{\otimes n}$.

The symmetric group $\Sym_n$ also acts on $V^{\otimes n}$ by place
permutations and it acts on $I(d;n)$ by permuting components;
indeed, $e_\a w=e_{\a w}$ for $\a\in I(d;n)$ and $w\in\Sym_n$.
Jimbo showed how to deform the action of $\Sym_n$ to give an action of
$\H_q(\Sym_n)$ on $V^{\otimes n}$. 

Recall that $\Sym_n$ is generated by $t_1,\dots,t_{n-1}$, where
$t_i=(i,i+1)$ for $i=1,\dots,n-1$.  Let $\Lambda^+(d;n)$ be the set of
partitions in $\Lambda(d;n)$.

\begin{Theorem}
[Jimbo~\cite{Jimbo}]\label{Jimbo}
Assume that $\H_q(\Sym_n)$ and $U_v(\gl_d)$ are defined over $\Q(v)$.
\begin{enumerate}
\item There is a unique $\H_q(\Sym_n)$-module structure on 
$V^{\otimes n}$ such that
$$e_\a\tilde T_j=\begin{cases} ve_\a,&\If a_j=a_{j+1},\\
                  \phantom ve_{\a t_j},&\If a_j>a_{j+1},\\
                  \phantom ve_{\a t_j}+(v-v^{-1})e_\a,&\If a_j<a_{j+1}
\end{cases}$$
for $j=1,\dots,n-1$ and $\a\in I(d;n)$.
\item The algebras $\H_q(\Sym_n)$ and
$\rho_n\(U_v(\gl_d)\)$ are mutually the full centralizer
algebras for each other for their actions on $V^{\otimes n}$.
Moreover, 
$$V^{\otimes n}\cong \bigoplus_{\lambda\in\Lambda^+(d;n)}
        W^\lambda\otimes S^\lambda$$
as an $\(\Schur_q(d;n),\H_q(\Sym_n)\)$-bimodule.
\end{enumerate}\end{Theorem}

It is not hard to see that there is an isomorphism 
$V^{\otimes n}\cong\oplus_{\lambda\in\Lambda(d;n)}M^\lambda$ of
$\H$-modules; consequently, $\rho_n\(U_v(\gl_d)\)\cong \Schur_q(d;n)$.
In part (ii), $W^\lambda$ is a Weyl module for $U_v(\gl_d)$; by what
we have just said this is the same as a Weyl module for the $q$-Schur
algebra $\Schur(d;n)$.

Actually, this is a slight modification of Jimbo's original action of
$\H_{v^2}(\Sym_n)$ on $V^{\otimes n}$; this action comes from
Du-Parshall-Wang~\cite{DPW}.

The proof of Theorem~\ref{Jimbo} is straightforward. Checking the
relations it is easy to see that $V^{\otimes n}$ is an $\H$-module
and that the actions of $\H$ and $U_v(\gl_d)$ commute. The double
centralizer property can be proved using a highest weight argument to
decompose $V^{\otimes n}$ as a $U_v(\gl_d)$-module.

Notice that Theorem~\ref{Jimbo} is stated over the rational function
field $\Q(v)$. Using the work of Be\uacc{\i}linson, Lusztig and
MacPherson \cite{BLM}, Du~\cite{Du:WeylRecip} showed that when $d\ge
n$ Theorem~\ref{Jimbo} holds over the Laurent polynomial ring
$\A=\Z[v,v^{-1}]$, where we replace $U_v(\gl_d)$ with its Lusztig
$\A$-form $U_{\A}(\gl_d)$; for the case $d<n$ see
\cite{DPS:qWeylRecip}. Further, if $d\ge n$ then
$\H_q(\Sym_n)\cong\End_{U_{\A}(\gl_d)}(V^{\otimes n})$.

We remark that Doty and Giaquinto~\cite{DG:pres} have recently used
the surjection $U(\gl_d)\To\Schur_q(d;n)$ to give a presentation of
the $q$-Schur algebras over $\Q(v)$; see also \cite{DP:pres}. No such
presentation is known for the cyclotomic $q$-Schur algebras.

Now we indicate how Sakamoto and Shoji~\cite{SakSho} have generalized
Theorem~\ref{Jimbo} to the cyclotomic case. We extend the coefficient
ring for all our algebras to the rational function field
$\Q(v,Q_1,\dots,Q_r)$, where $v,Q_1,\dots,Q_r$ are indeterminates.

Fix positive integers $d_1,\dots,d_r$ with $d=d_1+\dots+d_r$
and let $\gamma\map{\{1,\dots,d\}}\{1,\dots,r\}$ be the map
such that $\gamma(a)=s$ if $s$ is minimal such that 
$a\le d_1+\dots+d_s$ and let $V_s$ be the subspace of $V$ with basis
$\set{e_a|\gamma(a)=s}$, for $1\le s\le r$, and let
$\g=\gl_{d_1}(V_1)\oplus\dots\oplus\gl_{d_r}(V_r)$. We consider
$U_v(\g)$ as a Levi subalgebra of $U_v(\gl_d)$ in the natural way.
Then $V^{\otimes n}$ is a $U_v(\g)$-module by restriction; let
$\rho_{n,r}\map{U_v(\g)}\End(V^{\otimes n})$ be the corresponding
representation of $U_v(\g)$.

In order to extend the action of $\H_q(\Sym_n)$ on $V^{\otimes n}$ to
an action of $\H$ define linear operators $\varpi$ and $S_j$ on
$V^{\otimes n}$ by
$$e_\a\varpi=Q_{\gamma(a_1)}e_\a\quad\And\quad
  e_\a S_j=\begin{cases} e_\a\tilde T_j,&\If \gamma(a_{j-1})=\gamma(a_j),\\
                 e_{\a t_j},&\Otherwise,
\end{cases}$$
for $\a\in I(d;n)$ and $1\le j<n$.

Let $\Lambda(d_1,\dots,d_r;n)$ be the set of multicompositions
$\lambda$ of $n$ such that $|d_s|=|\lambda^{(s)}|$ for
$1\le s\le r$ and let $\Lambda^+(d_1,\dots,d_r;n)$ be the set of
multipartitions in~$\Lambda(d_1,\dots,d_r;n)$. We warn the reader that
$\Lambda^+(d_1,\dots,d_r;n)\ne\(\Lambda(d_1,\dots,d_r;n)\)^+$,
in the sense of Definition~\ref{Schur defn} --- unless $d_s\ge n$ for $1\le s<r$.

The irreducible representations of $U_v(\g)$ can be parametrized by
multipartitions in $\Lambda^+(d_1,\dots,d_r;n)$ for $n\ge0$. If
$\lambda\in\Lambda^+(d_1,\dots,d_r;n)$ let $W(\lambda)$ be the
corresponding Weyl module for~$U_v(\g)$. If $r=1$ then
$W(\lambda)\cong W^{\lambda}$ as $\Schur_q(d;n)$-modules; however, in
general, $W^\lambda$ and $W(\lambda)$ are not isomorphic even as
vector spaces.

\begin{Theorem}
[Sakamoto and Shoji~\cite{SakSho}]\label{SakSho}
Assume that $\H$ and $U_v(\g)$ are defined over the field
$\Q(v,Q_1,\dots,Q_r)$.
\begin{enumerate}
\item The action of $\H_q(\Sym_n)$ on $V^{\otimes n}$ extends to an
action of $\H$ on $V^{\otimes n}$ via
$$e_\a T_0=e_\a\varpi S_1\dots S_{n-1}
          \tilde T_{n-1}^{-1}\dots\tilde T_1^{-1}$$
for all $\a\in I(d;n)$.
\item The algebras $\H$ and
$\Schur^{(r)}_n\cong\rho_{n,r}\(U_v(\g)\)$ are
mutually the full centralizer algebras for the others
action on $V^{\otimes n}$. Moreover, 
$$V^{\otimes n}\cong \bigoplus_{\lambda\in\Lambda^+(d_1,\dots,d_r;n)}
        W(\lambda)\otimes S^\lambda$$
as an $\(\Schur^{(r)}_n,\H\)$-bimodule.
\end{enumerate}\end{Theorem}

Sakamoto and Shoji were guided in part by Ariki, Terasoma and
Yamada~\cite{ATY} who considered the special case when
$d_1=\dots=d_r=1$. The proof of part (i) of the theorem is a long
calculation building on Theorem~\ref{Jimbo}(i). Once again, part (ii)
is a highest weight computation.

Sakamoto and Shoji also note that part (i) of Theorem~\ref{SakSho} is
true over an arbitrary integral domain. Using this observation they
gave another proof that $\H$ is free of rank $|W_{r,n}|$
(Theorem~\ref{AK basis}).

Ariki~\cite{Ariki:Schur} asked whether Theorem~\ref{SakSho}(ii) is
true over an arbitrary field; he was particularly interested in
knowing
when the dimension of $\rho_{r,n}\(U_v(\g)\)=\End_\H(V^{\otimes n})$
is independent of $R$, $q$ and $\q$. Ariki found an example which
showed that in general the dimension of $\rho_{r,n}\(U_v(\g)\)$  does
depend upon these choices; nonetheless, he was able to prove the
result below.

Let $U_\A(\g)$ be the Kostant-Lusztig $\A$-form of $U_v(\g)$ and set 
$U_{R,v}(\g)=U_v(\g)\otimes_\A R$, where $R$ is an integral domain. We
also consider $V$ to be the free $R$-module with basis
$\{e_1,\dots,e_d\}$. Finally, define
$\Schur^{(r)}_n=\End_\H(V^{\otimes n})$.

\begin{Theorem}
[Ariki~\cite{Ariki:Schur}]\label{Ariki:Schur}
Suppose that $R$ is an integral domain and that $q=v^2,Q_1,\dots,Q_r$ are
elements of $R$ such that
$$P_n(q,\q)=\prod_{1\le i<j\le r}\prod_{-n<d<n}(q^dQ_i-Q_j)$$
is invertible in $R$. Then the following hold.
\begin{enumerate}
\item Suppose that $d_s\ge n$ for all $s$. Then there is an 
isomorphism of $R$-algebras
$$\Schur^{(r)}_n\cong
   \bigoplus_{\substack{n_1,\dots,n_r\ge0\\n_1+\dots+n_r=n}}
      \Schur_q(d_1;n_1)\otimes\dots\otimes\Schur_q(d_r;n_r).$$
In particular, $\Schur^{(r)}_n$ is free as an $R$-module and its rank
is independent of the choice of $R$ or the parameters 
$v,Q_1,\dots,Q_r$.
\item The algebra $\Schur^{(r)}_n$ is a quotient of $U_{R,v}(\g)$.
\item Assume that $d_s\ge n$ for $1\le s\le r$. Then
$\End_{U_{R,v}(\g)}(V^{\otimes n})$ is Morita equivalent to
the algebra
$$\bigoplus_{\substack{n_1,\dots,n_r\ge0\\n_1+\dots+n_r=n}}
      \H_{q}(\Sym_{n_1})\otimes\dots\otimes\H_q(\Sym_{n_r}).$$
\end{enumerate}\end{Theorem}

Observe that $P_n(q,\q)$ is equal to the polynomial
$P_n(q,\q_1,\dots,\q_r)$ of Theorem~\ref{Morita} (that is,
$|\q_\alpha|=1$ for all~$\alpha$). Assume that $d_s\ge n$ for all $s$.
Then, by part~(i) and Theorem~\ref{SMorita}, the algebra
$\Schur^{(r)}_n$ is Morita equivalent to the cyclotomic Schur algebra
$\Schur(\Lambda_{r,n})$.  Similarly, by part~(iii) and
Theorem~\ref{Morita}, if $d_s\ge n$ for all~$s$ then
$\End_{U_{R,v}}(V^{\otimes n})$ is Morita equivalent to~$\H$. Hence,
up to Morita equivalence, we have a complete analogue of Schur-Weyl
duality linking $U_v(\g)$, $\Schur(\Lambda_{r,n})$ and $\H$ in this
setting; however, note that this is really a type $A$ phenomenon and
is not genuinely `cyclotomic'.

Ariki also uses this result to compute the decomposition matrices of
the algebras $\Schur^{(r)}_n$ when $R=\Q$, $q\ne1$ and 
$P_n(q,\q)\ne0$. To do this he
uses part~(i) and the LT-conjecture~\cite{LT:Schur2} which gives an
extension of Theorem~\ref{Ariki}(ii) to the $q$-Schur algebras. The
LT-conjecture was proved by Varagnolo and Vasserot~\cite{VV}; 
see also Schiffmann~\cite{Schiff}.

Combining these results, the decomposition matrices of the cyclotomic
Schur algebras are known whenever $R$ is a field of characteristic
zero, $q\ne1$ and $P_n(q,\q)\ne0$. Actually, we do not need Ariki's
work to do this as we already obtain this result from
Theorem~\ref{SMorita} and~\cite{LT:Schur2,VV}. (Note that Ariki's
paper appeared before \cite{DM:Morita}.)

\subsection{Borel subalgebras}In this section we show that the
cyclotomic Schur algebras admit a ``triangular decomposition''. For
the Schur algebras this is a result of J.A.~Green~\cite{Green:Borel};
the cyclotomic case is due to Du and Rui~\cite{DuRui:Borel}. 

For simplicity we consider the case where $\Lambda=\Lambda_{r,n}$ is
the set of all multicompositions of $n$ of length at most $n$. Du and
Rui note that the general case can be deduced from this because if
$\Lambda$ is a saturated set of multicompositions then
$\Schur(\Lambda)$ is Morita equivalent to the subalgebra
$e\Schur(\Lambda_{r,n})e$ of $\Schur(\Lambda_{r,n})$, where $e$ is the
idempotent $\sum_{\lambda\in\Lambda^+}\phi_\lambda$.

Recall that $I(rn;n)=\set{(a_1,\dots,a_n)|1\le a_i\le rn}$. Then
$\Sym_n$ acts on $I(rn;n)$ by place permutations.  Given a
multicomposition $\lambda$ in $\Lambda_{r,n}$ let
$\bar\lambda=(\bar\lambda_1,\dots,\bar\lambda_{rn})$ be the
composition in $\Lambda(rn;n)$ with $\bar\lambda_i=\lambda^{(s)}_j$ 
if~$i=(s-1)n+j$. Define
$$\1_\lambda=(i_{\lambda,1},\dots,i_{\lambda,n})
            =(\underbrace{1,\dots,1}_{\bar\lambda_1},
              \underbrace{2,\dots,2}_{\bar\lambda_2}\dots,
              \underbrace{rn,\dots,rn}_{\bar\lambda_{rn}})\in I(rn;n).$$
Let $\succeq$ be the partial order on $I(rn;n)$ given by
$\a\succeq\mathbf b$ if $a_k\ge b_k$ for $1\le k\le rn$. Note that for
any $d\in\Sym_n$ if $\lambda$ and $\mu$ are multicompositions with
$\1_\lambda d\succeq\1_\mu$ then $\mu\gedom\lambda$.

Recall that for each multicomposition $\lambda\in\Lambda_{r,n}$ we
have a Young subgroup $\Sym_\lambda$ and that $\D_\lambda$ is the set of
minimal length coset right representatives for $\Sym_\lambda$ in
$\Sym_n$.  Moreover, if $\mu$ is another multicomposition then
$\D_{\lambda\mu}=\D_\lambda\cap\D_\mu^{-1}$ is a set of minimal length
$(\Sym_\lambda,\Sym_\mu)$-double coset representatives. For each
$d\in\D_{\lambda\mu}$ define $\phi^d_{\lambda\mu}$ to be the
$R$-linear endomorphism of $\bigoplus_\alpha M^\alpha$ determined by
$$\phi^d_{\lambda\mu}(m_\alpha h)=\delta_{\alpha\mu}
          \Big(\sum_{w\in\Sym_\lambda d\Sym_\mu} T_w\Big)
                u_\mu^+h$$
for all $\alpha\in\Lambda_{r,n}$ and all $h\in\H$. If 
$\1_\lambda d\succeq\1_\mu$ then 
$\phi^d_{\lambda\mu}\in\Schur(\Lambda_{r,n})$ by
\cite[Lemma~5.6]{DuRui:Borel}. In particular, if~$\nu\in\Lambda_{r,n}$
then $\phi^1_{\nu\nu}=\phi_\nu$ restricts to the identity map 
on $M^\nu$ (and is zero on $M^\alpha$ for~$\alpha\ne\nu$).

Finally, given multicompositions $\lambda$ and $\mu$ in $\Lambda_{r,n}$
let
$$\Omega_{\lambda\mu}
     =\set{d\in\D_{\lambda\mu}|\1_\lambda d\succeq\1_\mu}.$$
Define $\Schur^\pm(\Lambda_{r,n})$ to be the two $R$-submodules of
$\Schur(\Lambda_{r,n})$ spanned by 
$\set{\phi^d_{\lambda\mu}|d^\mp\in\Omega_{\lambda\mu}}$.
We can now state the main result.

\begin{Theorem}
[Du and Rui~\cite{DuRui:Borel}]\label{borel}
Suppose that $R$ is an integral domain ring.
\begin{enumerate}
\item The two $R$-modules $\Schur^\pm(\Lambda_{r,n})$ are 
subalgebras of $\Schur(\Lambda_{r,n})$.
\item $\Schur^\pm(\Lambda_{r,n})$ is free as an $R$-module with basis
$$\set{\phi^d_{\lambda\mu}|\lambda,\mu\in\Lambda_{r,n}
                          \And d^\mp\in\Omega_{\lambda\mu}}.$$
\item $\Schur(\Lambda_{r,n})$ has a triangular decomposition
$$\Schur(\Lambda_{r,n})
     =\Schur^-(\Lambda_{r,n})\cdot\Schur^+(\Lambda_{r,n})
     =\Schur^-(\Lambda_{r,n})\cdot
      \Big(\sum_{\nu\in\Lambda_{r,n}}R\phi_\nu\Big)\cdot
      \Schur^+(\Lambda_{r,n}).$$
Thus, $\set{\phi^d_{\lambda\mu}\phi^e_{\mu\nu}|
    \lambda,\mu,\nu\in\Lambda_{r,n}, d\in\D_{\lambda\mu}\And
                               e^{-1}\in\D_{\mu\nu}}$
is a basis of $\Schur(\Lambda_{r,n})$.
\end{enumerate}\end{Theorem}        

Du and Rui call $\Schur^-(\Lambda_{r,n})$ and 
$\Schur^+(\Lambda_{r,n})$ the {\sf Borel subalgebras} 
of~$\Schur(\Lambda_{r,n})$. Surprisingly, the Borel subalgebras of the
cyclotomic Schur algebras are isomorphic to the Borel subalgebras of
the $q$-Schur algebras; hence, they are really type $A$ algebras.

The right hand side of part (iii) is
written so as to suggest the triangular decomposition of quantum
groups; however, this is slightly misleading because 
$\phi^d_{\lambda\mu}\(\sum_\nu r_\nu\phi_\nu\)\phi^e_{\sigma\tau}
  =\delta_{\mu\sigma}r_\mu\phi^d_{\lambda\mu}\phi^e_{\sigma\tau},$
for $r_\nu\in R$.

Du and Rui are able to say quite a lot about the representation
theory of these subalgebras. Because $\Schur^\pm(\Lambda_{r,n})$ are
quasi-hereditary, they have standard modules and costandard modules;
denote these by $\Delta^\pm(\mu)$ and $\nabla^\pm(\mu)$ respectively,
for $\mu\in\Lambda_{r,n}$. Also, if $\mu\in\Lambda_{r,n}^+$ then the
Weyl module $W^\mu=\Delta(\mu)$ is a standard module of
$\Schur(\Lambda_{r,n})$ and its contragredient dual
$(W^\mu)^*=\nabla(\mu)$ is a costandard module (duality with respect
to~$*$).

\begin{Theorem}
[Du and Rui]\label{induced modules}
Suppose that $R$ is a field.
\begin{enumerate}
\item The Borel subalgebras $\Schur^-(\Lambda_{r,n})$ and
$\Schur^+(\Lambda_{r,n})$ are quasi-hereditary, with respect to the
poset $\Lambda_{r,n}$. Moreover, $\Schur^-(\Lambda_{r,n})$ and
$\Schur^+(\Lambda_{r,n})$ are Ringel dual to each other.
\item
\begin{enumerate}
\item Each costandard module of $\Schur^-(\Lambda)$ is one
  dimensional and, hence, simple; moreover, every simple module appears
  this way.
  \item Dually, each standard module of $\Schur^-(\Lambda_{r,n})$ is a
  projective indecomposable $\Schur^-(\Lambda_{r,n})$--module.
  \item Explicitly, if $\mu\in\Lambda_{r,n}^+$ then
  $\Delta^-(\mu)=\Schur^-(\Lambda_{r,n})\phi_\mu$ and
  $\nabla^-(\mu)=\Delta^-(\mu)/\Rad\Delta^-(\mu)$; moreover,
  $\set{\phi_\mu|\mu\in\Lambda_{r,n}}$ is a complete
  set of primitive idempotents in $\Schur^-(\Lambda_{r,n})$.
\end{enumerate}
\item Suppose that $\mu\in\Lambda_{r,n}$. Then
\begin{align*}
\Schur(\Lambda_{r,n})\otimes_{\Schur^+(\Lambda_{r,n})}\Delta^+(\mu)
     &\cong\begin{cases} \Delta(\mu),&\If\mu\in\Lambda_{r,n}^+,\\
                     0,&\Otherwise,
\end{cases}
\intertext{and}
\Hom_{\Schur^-(\Lambda_{r,n})}\(\Schur(\Lambda_{r,n}),\nabla^-(\mu)\)
     &\cong\begin{cases} \nabla(\mu),&\If\mu\in\Lambda_{r,n}^+,\\
                     0,&\Otherwise,
\end{cases}\end{align*}
\end{enumerate}\end{Theorem}

Ringel duality interchanges the standard and costandard modules of
$\Schur^-(\Lambda_{r,n})$ and $\Schur^+(\Lambda_{r,n})$, so part (ii)
also describes the simple and projective
$\Schur^+(\Lambda_{r,n})$-modules.

Du and Rui also give the dimensions of the standard and costandard
modules for the Borel subalgebras and show that the Borel subalgebras are
full tilting modules for the Ringel duality. 

\subsection{Tilting modules}

Let $A$ be a quasi-hereditary algebra (see~\cite{CPS:qh,Donkin:book}),
and let $\Lambda^+$ be its poset of weights. Then for each
$\lambda\in\Lambda^+$ we have a standard module $\Delta(\lambda)$, a
costandard module $\nabla(\Lambda^+)$ and a simple module
$L(\lambda)$. The simple module $L(\lambda)$ is the head of
$\Delta(\lambda)$ and the simple socle of $\nabla(\lambda)$; further,
$\nabla(\lambda)$ is the contragredient dual of $\Delta(\lambda)$ if
$A$ possesses a suitable involution.

Let $\F(\Delta)$ be the full subcategory of $A\mod$ consisting of those
modules which have a $\Delta$-filtrations; thus $X\in\F(\Delta)$ if $X$
has a filtration $X=X_1\supset X_2\supset\dots\supset X_m\supset0$
with $X_i/X_{i+1}\cong\Delta(\lambda_i)$ for $1\le i\le m$. If
$X\in\F(\Delta)$ and $\lambda\in\Lambda^+$ let 
$[X{:}\Delta(\lambda)]
   =\#\set{1\le i\le m|X_i/X_{i+1}\cong\Delta(\lambda)}$;
this is independent of the choice of filtration because the
equivalence classes of standard modules are a basis of the Grothendieck
group of $A$. Similarly, let $\F(\nabla)$ be the full subcategory of
$A$-modules which have a $\nabla$-filtration. 

Ringel~\cite{Ringel} has proved that for each $\lambda\in\Lambda^+$
there is a unique indecomposable $A$-module
$T(\lambda)\in\F(\Delta)\cap\F(\nabla)$ such that
$[T(\lambda){:}\Delta(\lambda)]=1$ and
$[T(\lambda){:}\Delta(\mu)]\ne0$ only if $\mu\ge\lambda$; we call
$T(\lambda)$ a (partial) {\sf tilting module} for $A$. Moreover, every
module in $\F(\Delta)\cap\F(\nabla)$ is isomorphic to a direct sum of
tilting modules.

If $\Lambda$ is saturated then the cyclotomic Schur algebra
$\Schur(\Lambda)$ is quasi-hereditary by Theorem~\ref{quasi}, so we
may ask for a description of the tilting modules of $\Schur(\Lambda)$.
When $r=1$ Donkin~\cite{Donkin:tilt,Donkin:book} determined the
tilting modules of the $q$-Schur algebras. To describe this, recall
from the previous section that $\Schur_q(d;n)=\End_\H(V^{\otimes n})$.
Donkin showed that when $d\ge n$ the tilting modules of
$\Schur_q(d;n)$ are precisely the indecomposable direct summands of
the exterior powers $\wedge^\lambda V
=\wedge^{\lambda_1}V\otimes\dots\otimes\wedge^{\lambda_d}V$.  For a
different approach to the tilting modules of the $q$-Schur algebras
see \cite{DPS:qWeylRecip}.

Even though we do not know how to describe $\oplus_\mu M^\mu$ as a
tensor product the exterior powers of $\Schur(\Lambda)$ still admit a
similar description.  In introducing $M^\lambda$ we said that it
should be thought of as an induced trivial module; the analogue of an
induced sign representation for $\H$ is the module
$N^\lambda=n_\lambda\H$, where $n_\lambda=y_\lambda
u_\lambda^-=u_\lambda^-y_\lambda$ and 
$$y_\lambda=\sum_{w\in
\Sym_\lambda} (-q)^{-\len(w)}T_w\quad\And\quad
u_\lambda^-=\prod_{s=1}^{r-1}
\prod_{k=1}^{|\lambda^{(s+1)}|+\dots+|\lambda^{(r)}|}(L_k-Q_s).$$ For
each multipartition $\lambda$ let 
$E^\lambda=\Hom_\H\(\bigoplus_{\mu\in\Lambda}M^\mu,N^\lambda\)$.  Then
$E^\lambda$ is a right $\Schur(\Lambda)$-module and we have the
following.

\begin{Theorem}
[Mathas~\cite{M:tilting}]
Suppose that $R$ is a field, $Q_s\ne0$, for $1\le s\le r$, and
that~$\Lambda$ is a saturated set of multicompositions containing
$\omega$.  Then the tilting modules of~$\Schur(\Lambda_{r,n})$ are the
indecomposable summands of the modules
$\set{E^\lambda|\lambda\in\Lambda^+}$.
\label{tilting}\end{Theorem}

The requirement that $Q_s\ne0$, for all $s$, essentially comes from
Proposition~\ref{symmetric}(ii); this assumption is needed to show that
the Specht modules and dual Specht modules are isomorphic. In turn,
this is used to show that $E^\lambda$ is self--dual. The key tool in
the proof of Theorem~\ref{tilting} is the use of Specht filtrations
and dual Specht filtrations of $\H$--modules; this is a bit surprising
because Specht filtrations are generally not as good as Weyl
filtrations (since $D^\lambda$ can be zero).

The tilting modules of $\Schur(\Lambda)$ have all of the expected
properties. For example,
$[T(\lambda){:}\nabla(\mu)]=[\Delta(\mu'){:}L(\lambda')]$ for all
$\lambda,\mu\in\Lambda^+$. (Here $\mu'$ is the multipartition
conjugate to $\mu$.) Furthermore, the Ringel dual of~$\Schur(\Lambda)$
is the algebra
$\Schur'(\Lambda)=\End_\H\(\bigoplus_{\mu\in\Lambda}N^\mu\)$ and
$\Schur'(\Lambda)\cong \Schur(\Lambda)$.

The theory of Young modules for $\H$ (cf.~\cite{James:trivsour}), is
also developed in \cite{M:tilting}. The Young modules are the
indecomposable direct summands of the modules $M^\lambda$ and the
$N^\lambda$, for $\lambda$ a multicomposition of $n$; they are indexed
by the multipartitions of $n$. The Young modules are the common image
under the Schur functor of the corresponding indecomposable
projective, injective or tilting module for either of the algebras
$\Schur(\Lambda)$ or $\Schur'(\Lambda)$.  


\section{Some open problems}

In this final chapter we discuss some open problems for the
Ariki-Koike algebras and the cyclotomic Schur algebras. We are mostly
interested in the connections between the representation theory of the
Ariki-Koike algebras and cyclotomic Schur algebras with the
representation theory of the finite groups of Lie type. 

\begin{Problem} Prove the conjectures of Brou\'e and Malle~\cite{BM:cyc}.
\end{Problem}

We also pose the more general (and more vague) problem.

\begin{Problem} Find a link between the representation theory of the
cyclotomic Schur algebras and the modular representation theory of the
finite groups of Lie type.
\end{Problem}

At best, there is only circumstantial evidence for such a connection
when~$r>2$. If we believe in the conjectures of the Brou\'e school
then there are strong ties between the representation theory of
cuspidal representations of $\GL_d(q)$ in characteristic zero, so it
is not unreasonable to expect that the modular theory of the
cyclotomic Schur algebras and Ariki-Koike algebras also carry
information about the modular representations of $\GL_d(q)$.

\subsection{Quantum groups and geometry}

The results of Ariki and Sakamoto and Shoji from section~5.4 show
that in some circumstances the module categories of the Ariki-Koike
algebras and the cyclotomic Schur algebras are connected with the
module categories of Levi subalgebras of $U_v(\gl_d)$. Unfortunately,
these results apply only in cases where the Ariki-Koike algebras are
Morita equivalent to direct sums of tensor products of Iwahori-Hecke
algebras of type $A$ and when the cyclotomic $q$-Schur algebras were
Morita equivalent to direct sums of tensor products of $q$-Schur
algebras.

\begin{Problem} Realize the cyclotomic Schur algebras as a quotient of a
quantum group $U_\A(\g)$ over an arbitrary integral domain.
\end{Problem}

We could ask for a generalization of the results of Sakamoto and
Shoji~(Theorem~\ref{SakSho}) and Ariki~(Theorem~\ref{Ariki:Schur});
however, as the conjectures of Brou\'e's school only ask for a derived
equivalence it seems to me that we cannot expect something so simple
here.

That the cyclotomic Schur algebras might be realizable as a quotient
of a quantum group is suggested by the cyclotomic Jantzen sum formula
(Theorem~\ref{sum formula}) and by the existence of the Borel
subalgebras and the triangular decomposition of $\Schur(\Lambda)$
(Theorem~\ref{borel}). Both of these results hint at connections with
quantum groups and at some undiscovered geometry. 

Note also that the existence of the Borel subalgebras allows us to
consider the dual Weyl modules of the cyclotomic Schur algebras as
induced modules and so gives us cohomological techniques to play with.

\subsection{Tensor products}

First consider the case $r=1$. If $\lambda$ is a partition of $n$ and
$\mu$ is a partition of $m$ then $S^\lambda\otimes S^\mu$ is a module for
the Hecke algebra $\H_q(\Sym_n)\otimes\H_q(\Sym_m)$. We can identify
$\H_q(\Sym_n)\otimes\H_q(\Sym_m)$ with the subalgebra
$\H_q(\Sym_{(n,m)})$ of $\H_q(\Sym_{n+m})$ where
$\Sym_{(n,m)}=\Sym_n\times\Sym_m$.  Thus, $\H_q(\Sym_{n+m})$ is a free
$\H_q(\Sym_n)\otimes\H_q(\Sym_m)$-module and we can define the
$\H_q(\Sym_{n+m})$-module
$$S^\lambda\otimes S^\mu=(S^\lambda\boxtimes S^\mu)
         \otimes_{\H_q(\Sym_n)\otimes\H_q(\Sym_m)}\H_q(\Sym_{n+m}).$$
When $\H_q(\Sym_n)$ is semisimple, this decomposes as a direct sum of
Specht modules according to the Littlewood-Richardson rule.

In the case of the $q$-Schur algebras it is even easier. If $\lambda$
and $\mu$ are both partitions of length at most $d$ then the Weyl
modules $W^\lambda$ and $W^\mu$ are homogeneous polynomial
representations for $U_v(\gl_d)$ of degree $n$ and $m$ respectively;
therefore, $W^\lambda\otimes W^\mu$ is a polynomial representation of
$U_v(\gl_d)$ of degree $n+m$ --- since $U_v(\gl_d)$ is a Hopf algebra.
Hence, $W^\lambda\otimes W^\mu$ is an $\Schur_q(d;n+m)$-module since
$\Schur(d;N)\mod$ is the category of polynomial representations of
$U_v(\gl_d)$ of homogeneous degree $N$. Again, in the semisimple case
the decomposition of $W^\lambda\otimes W^\mu$ into irreducibles is
given by the Littlewood-Richardson rule.

When we try and extend either of these constructions to the cyclotomic
case we run into problems. First, for the Ariki-Koike algebras there
is no obvious way to consider
$\H_{q,\q}(W_{r,n})\otimes\H_{q,\q'}(W_{s,m})$ as a free submodule of
$\H_{q,\q\cup\q'}(W_{t,n+m})$ for any $t$, unless $rs=0$.  Secondly,
we do not have an interpretation of the module category of a
cyclotomic Schur algebra in terms of homogeneous representations of a
quantum group.

\begin{Problem} Find a good tensor product operation for the categories
$\H\mod$ and $\Schur(\Lambda)\mod$.
\end{Problem}

Of course, a strong enough link with quantum groups would give us this
for free.  The correct approach is probably via the affine Hecke
algebra (or possibly the work of Shoji~\cite{Sho:Frob}).

If we knew how to take tensor products of modules for the cyclotomic
Schur algebras then we could try and solve the following problem.

\begin{Problem} Find an analogue of the Steinberg tensor product
theorem for the cyclotomic Schur algebras.
\end{Problem}

Evidence for the existence of such a result, as well as an indication
of what it might look like, are given by Uglov's~\cite{Uglov} action of
the Heisenberg algebra upon the generalized Fock spaces.


\subsection{Decomposition numbers at roots of unity} The decomposition
numbers of the Ariki-Koike algebras are known in characteristic zero,
thanks to Ariki's theorem and the work of Uglov (assuming that
$Q_s\ne0$ for any $s$); see Corollary~\ref{decomp nos}. 

\begin{Problem} Compute the decomposition numbers of the cyclotomic
$q$-Schur algebras in characteristic zero.
\end{Problem}

By Theorem~\ref{Ariki} the decomposition matrix of $\H_q(\Sym_n)$ can be
computed from the canonical basis of $L_v(\Lambda_0)$. The easiest way
to compute the canonical basis of $L_v(\Lambda_0)$ is to work in the
Fock space $\F$, an infinite rank free $\C[v,v^{-1}]$-module with
a basis given by the set of all partitions of all integers.  Leclerc
and Thibon's idea~\cite{LT:Schur2} was to define a canonical basis on
the whole of the Fock space; they did this using the action of a
Heisenberg algebra on~$\F$. Leclerc and Thibon conjectured
that the decomposition matrices of the $q$-Schur algebra were given
by computing the canonical basis of $\F$ and then specializing
$v=1$; this was proved by Varagnolo and Vasserot~\cite{VV}.

Hence, this problem has been solved when $r=1$. Furthermore, as
remarked in section~5.4, when $P_n(q,\q)\ne0$ we also know the answer
because by Theorem~\ref{SMorita} $\Schur(\Lambda)$ is Morita
equivalent to a direct sum of tensor products of $q$-Schur algebras.

Now, the decomposition matrices of the Ariki-Koike algebras in
characteristic zero are obtained by computing the canonical basis of
highest weight modules $L_v(\Lambda)$, for the various dominant
weights $\Lambda$. This time, $L_v(\Lambda)$ embeds in a generalized
Fock space $\F_\Lambda$ and Uglov has shown how to compute a canonical
basis for the whole of this space; this gives a canonical basis
element for each multipartition. For $n\ge0$ the canonical basis of
$\F_\Lambda$ at $v=1$ gives a square unitriangular matrix, indexed by
the multipartitions of $n$, which contains the decomposition matrix of
the Ariki-Koike algebra $\H_n$ as a submatrix (delete those columns
corresponding to the multipartitions~$\lambda$ with $D^\lambda=0$);
compare with Corollary~\ref{submatrix}. The indexing of the rows and
columns is wrong; however, modulo the difference in labeling this
should be the decomposition matrix of $\Schur(\Lambda_{r,n})$.

\subsection{Dipper-James theory}Let $q$ be a prime power and let
$\GL_n(q)$ be the general linear group over a field with $q$ elements.
Dipper and James~\cite{DJ:Schur} proved that the decomposition matrix
of $\GL_n(q)$ in non-defining characteristic is completely determined
by the decomposition matrix of the $q^d$--Schur algebras, for $d\ge1$.
Recently Brundan, Dipper, and Kleshchev~\cite{BDK} have rewritten this
theory using cuspidal algebras. They also make the Dipper-James
result on decomposition matrices much more explicit;
see~\cite[Theorem~4.4d]{BDK}.

There is no analogue of this theory for the cyclotomic $q$-Schur
algebras.  The best results in this direction were obtained by Gruber
and Hiss~\cite{GrHi} who, for linear primes, worked with a Morita
equivalent version of the cyclotomic Schur algebras when $r=2$ (type
$B$), to give weaker results for other finite reductive groups
$G_n(q)$. See the survey article of Dipper, Geck, Hiss and
Malle~\cite{DGHM} for the current status of this theory.

\section*{Acknowledgements} I would like to thank Michel Brou\'e for
explaining his conjectures on cyclotomic Hecke algebras to me, Gunter
Malle and Jean Michel for many useful comments and corrections and Jie
Du for some further references. Part of this article was written at
the University of Leicester; I thank them, and in particular Steffen
K\"onig, for their hospitality. Finally, I would like to thank the
organizers for such a good meeting and for giving me the chance to
speak.


\let\u\uacc   

\end{document}